\documentclass[final]{siamart190516}
\usepackage[utf8]{inputenc}
\usepackage[T1]{fontenc}
\usepackage{graphicx}
\usepackage[english]{babel}
\usepackage{amsmath,amssymb,amsbsy}
\usepackage{float}
\usepackage{bm}
\usepackage{multirow}
\usepackage{url}
\usepackage{hyperref}
\usepackage{algorithmic}

\newcommand{\R}{\mathbb{R}}
\newcommand{\N}{\mathbb{N}}
\newcommand{\T}{\mathbb{T}}
\newcommand{\Z}{\mathbb{Z}}
\renewcommand{\P}{\mathbb{P}}
\renewcommand{\L}{\mathbb{L}}

\renewcommand{\d}{\mathrm{d}}

\newsiamremark{remark}{Remark}

\Crefname{ALC@unique}{Line}{Lines} 

\headers{Adaptive sparse trigonometric interpolation}{Zack Morrow and Miroslav Stoyanov}

\title{A method for dimensionally adaptive sparse trigonometric interpolation of periodic functions}

\author{Zack Morrow\thanks{Department of Mathematics, North Carolina State University, Raleigh, NC 27695 (\href{mailto:zbmorrow@ncsu.edu}{zbmorrow@
ncsu.edu}).} 
\and Miroslav Stoyanov\thanks{Corresponding author. Computational and Applied Mathematics Group, Oak Ridge National Laboratory, Oak Ridge, TN 37830 (\email{stoyanovmk@ornl.gov}).}
}

\begin{document}
\maketitle 

\begin{abstract}
	We present a method for dimensionally adaptive sparse trigonometric interpolation of multidimensional periodic functions belonging to a smoothness class of finite order. This method targets applications where periodicity must be preserved and the precise anisotropy is not known {\em a priori}. To the authors' knowledge, this is the first instance of a dimensionally {\em adaptive} sparse interpolation algorithm that uses a {\em trigonometric interpolation} basis. The motivating application behind this work is the adaptive approximation of a multi-input model for a molecular potential energy surface (PES) where each input represents an angle of rotation. Our method is based on an anisotropic quasi-optimal estimate for the decay rate of the Fourier coefficients of the model;  a least-squares fit to the coefficients of the interpolant is used to estimate the anisotropy. Thus, our adaptive approximation strategy begins with a coarse isotropic interpolant, which is gradually refined using the estimated anisotropic rates. The procedure takes several iterations where ever-more accurate interpolants are used to generate ever-improving anisotropy rates. We present several numerical examples of our algorithm where the adaptive procedure successfully recovers the theoretical ``best'' convergence rate, including an application to a periodic PES approximation. An open-source implementation of our algorithm resides in the Tasmanian UQ library developed at Oak Ridge National Laboratory.
\end{abstract} 



\section{Introduction}
\label{sec:intro}

Consider the approximation of a periodic multidimensional function (i.e., model) $f : \T^d \to \R$, where $\T = [0,1]$ represents the unit interval and $f$ is both differentiable and periodic up to a given finite order. Since the domain in our motivating application is a torus, we chose the letter $\T$ to help reinforce the notion of periodicity. The order of differentiability and periodicity can vary with each dimension. The main challenge of approximating a computationally expensive model in a multidimensional (multi-input) setting is the rapid growth of required model simulations per dimension $d$, a phenomenon known as the \textit{curse of dimensionality}~\cite{Bellman_1961}.

Various techniques exist to mitigate, or in rare cases eliminate, the curse of dimensionality. Global and derivative-based sensitivity analysis enables the identification of non-influential parameters, thereby reducing the effective dimensionality of the problem to include only inputs and directions that contribute towards the model output variability~\cite{Hart_Alexanderian_Gremaud_2017,Stoyanov_Webster_2015}. Other methods seek to reduce the complexity of the target function $f$ by approximating it with functions that are in some sense ``simpler,'' e.g., \cite{Bungartz_Griebel_2004, EldredBurkardt2009, EldredWebsterConstantine2008, jakeman2011characterization, jakeman2012local, nobile2008sparse, pfluger2010spatially, Stoyanov_Webster_2016}. This could be done, for instance, by projection onto or interpolation within a polynomial or trigonometric function space, both of which use samples, i.e., the values of the target function for a set of independent inputs. Sampling methods are attractive because those can be wrapped around existing third-party or black-box models in a non-intrusive way, i.e., without the need to modify the original solver.

Let $\{\phi_{\bm \nu}\}_{\bm \nu \in \N^d}$ be an orthonormal basis for the Hilbert space where $f$ resides, e.g., $L^2(\T^d)$, and let $\Lambda \subset \N^d$ be finite.\footnote{This paper adopts the convention $\N = \{0,1,2,\dots\}$.} It is well known that orthogonal projection of $f$ onto $S = \text{span} \{\phi_{\bm \nu}\}_{\bm \nu \in \Lambda}$ yields the optimal $L^2$ error~\cite[p.~352]{Kreyszig_1978}; that is,
\begin{equation} 
c_{\bm \nu} = \langle f, \phi_{\bm \nu} \rangle_{L^2(\T^d)}
\quad \Rightarrow \quad 
\left\| f - \sum_{\bm \nu \in \Lambda} c_{\bm \nu} \phi_{\bm \nu} \right\|_{L^2(\T^d)} = \min_{g \in S} \, \left\| f - g \right\|_{L^2(\T^d)}.
\label{equ:c_nu}
\end{equation}
Here, $\{ c_{\bm \nu} \}_{\bm \nu \in \Lambda}$ are the optimal expansion coefficients, and $\langle \cdot, \cdot \rangle$ and $\| \cdot \|$ denote the $L^2(\T^d)$ inner product and norm respectively. In general, the integral coefficients in \cref{equ:c_nu} must be evaluated numerically, e.g., with a multidimensional numerical quadrature. Thus, projection methods often come at a high computational cost due to the large number of function samples necessary to approximate $c_{\bm \nu}$ to a sufficient accuracy which could far exceed the number of basis functions~\cite{Beck_et_al_2014,Stoyanov_Webster_2016,Todor_Schwab_2007}. In contrast, interpolation methods require a single sample per basis function, although the resulting approximation is not Hilbert-optimal. The interpolation error is bounded by the best projection error multiplied by a penalty term called the \textit{Lebesgue constant}, but the degradation in accuracy is usually offset by the reduction in computational cost. In particular, sparse interpolation methods~\cite{Smolyak_1963}, which are of focus of this paper, often have better overall convergence rate with respect to the number of samples~\cite{Stoyanov_Webster_2016}. 

Aiming to further improve the convergence rate of sparse-grid interpolation, many methods gauge the approximation error in order to determine the most important directions and spatial locations in which to sample next. Such procedures are known as \textit{adaptive refinement}, and the overall goal is to select the samples that would result in the fastest convergence rate. Bungartz and Griebel formulated this procedure as a knapsack problem in which they maximize the added accuracy subject to cost constraints at each refinement iteration until the interpolation error reaches a desired accuracy~\cite{Bungartz_Griebel_2004}. The greedy construction of an approximate knapsack set has demonstrated good performance in many applications; however, the exact solution to the knapsack problem remain intractable and the greedy approach is susceptible to stagnation and premature termination due to non-monotonic behavior of $c_{\bm \nu}$ when the grid is still fairly coarse (also called the {\em preasymptotic} regime).

Recent developments in quasi-optimal approximation utilize theoretical upper bounds on the decay rates of the expansion coefficients $c_{\bm \nu}$ (rather than the exact $c_{\bm \nu}$ themselves) in order to refine the approximation iteratively. For a fixed number of terms $M$ in the expansion, such approaches often result in a tighter error bound. Much of the previous work on quasi-optimal approximation is done in the context of Legendre expansions of holomorphic functions. Authors have deployed quasi-optimal approximation in the context of projection~\cite{Tran_Webster_Zhang_2017} and sparse-grid interpolation~\cite{Nobile_Tamellini_Tempone_2016, Stoyanov_Webster_2016}.

In contrast to previous work, this paper considers periodic functions, i.e., $f \in L^2(\mathbb{T}^d)$, such that, for each $1 \leq k \leq d$,
$$
\left\| \frac{\partial^m f}{\partial x_k^m} \right\|_{L^2(\T^d)} < \infty, \hspace{2em} \frac{\partial^m f(0)}{\partial x_k^m} = \frac{\partial^m f(1)}{\partial x_k^m}, \hspace{2em} \forall m \in \{0,1,\dots, n_k \}
$$
where $n_k$ is the order of periodicity in dimension $k$. The periodic boundary conditions lead to the natural choice of basis of trigonometric polynomials~\cite{Jackson_1930}. In this paper, we present a multidimensional sparse-grid trigonometric-basis adaptive interpolation technique for periodic functions with different degrees of smoothness in each direction. The smoothness affects the convergence rate, so we will use anisotropic grids with rates of anisotropy estimated using a least-squares fitting of the Fourier coefficients of the interpolant, which is similar to previous work in the context of total-degree, polynomial-based interpolation for holomorphic functions~\cite{Stoyanov_Webster_2016}. To the authors' knowledge, this is the first instance of a method for adaptive sparse-grid interpolation with a trigonometric basis.

The rest of the paper is organized as follows. In~\cref{sec:quasi_opt_general}, we derive the trigonometric quasi-optimal approximation space using theoretical bounds on the decay rates of the Fourier coefficients. In~\cref{sec:sparse_grids}, we describe anisotropic sparse trigonometric interpolation and present our adaptive refinement algorithm. In~\cref{sec:numerical_results}, we provide several numerical examples using both simple polynomials with known degrees of periodicity and the PES model. Finally,~\cref{sec:conclusion} offers a brief summary of our results.

\section{Quasi-optimal function space}
\label{sec:quasi_opt_general}
First, we consider the space of multidimensional periodic functions. Using upper bounds on the Fourier coefficients, we derive the quasi-optimal approximation space in the context of projection. From projection, we proceed to interpolation and derive the quasi-optimal interpolation space. We conclude by discussing how to estimate the anisotropic coefficients of the target function {\em on-the-fly}.

Let $H^n(\T) \subset C^n(\T)$, with $n \geq 0$, denote the space of $n$-times continuously differentiable functions $f:\T \to \R$ such that $f$ has $n$ periodic derivatives and $f^{(n+1)}$ is piecewise continuous with only finitely many jump discontinuities. Functions arising from science and engineering applications often satisfy this mild requirement of piecewise differentiability (see e.g., ~\cite{Levine_2014,Nance_Jakubikova_Kelley_2014} and~\cref{sec:numerical_results}). In the $d$-dimensional case, for $\bm n = (n_1, n_2, \cdots, n_d)$, we define
$$
H^{\bm n}(\T^d) = H^{n_1}(\T) \otimes \cdots \otimes H^{n_d}(\T)
$$
so that for any $f \in H^{\bm n}(\T^d)$, $1 \leq k \leq d$, and $(x_1, x_2, \cdots, x_d) \in \T^d$
$$
    f(x_1, \cdots, x_{k-1}, x, x_{k+1}, \cdots x_d) \in H^{n_k},
$$
i.e., restricting $f$ to a single dimension yields a function in $H^{n_k}(\T)$. Here, without loss of generality, we take the canonical $\T^d = [0, 1]^d$ since any arbitrary hypercube $\Gamma = \bigotimes_{k=1}^d [a_k,b_k]$ can be translated to $\T^d$ with a simple affine transformation.

The coefficients of the $L^2$-Fourier expansion of $f \in H^{\bm n}$ are defined as
\begin{equation}
c_{\bm j}(f) = \int_{\T^d} \exp(-2 \pi \text{i} \, \bm j \cdot \bm x) f(\bm x) \, \d \bm x,~~~~\bm j \in \Z^d,
\label{equ:fcoef_def}
\end{equation}
with $\text{i}^2 = -1$ and $\bm j \cdot \bm x = \sum_{k=1}^d j_k x_k$. In a single-dimensional context,
using Theorems 1.6, 4.4, and 4.5 from~\cite[pp.~4, 25]{Katznelson_2004} and trivial re-indexing, we obtain
\begin{equation}
|c_j(f)| \leq \frac{C(f)}{(1+|j|)^{n+2}}\, ,~~~~j \in \Z \, ,~f \in H^n(\T),
\label{equ:fcoef_bound_1d}
\end{equation}
for some constant $C(f)>0$ that depends on $f$. Furthermore, since $f^{(n+1)}$ has jump discontinuities, the bound in~\cref{equ:fcoef_bound_1d} is asymptotically sharp~\cite[p.~200]{Grafakos_2014}. In a multidimensional context, using the tensor-product structure of the space, we have 
\begin{equation}
|c_{\bm j}(f)| \leq \frac{C(f)}{\prod_{k=1}^d (1+|j_k|)^{n_k+2}},~~~~\bm j \in \Z^d,~f \in H^{\bm n}(\T^d) \, .
\label{equ:fcoef_bound_multidim}
\end{equation}

Function spaces like this have appeared in the literature as weighted Sobolev and Korobov spaces. In an early work on sparse trigonometric interpolation, Hallatschek~\cite{Hallatschek_1992} considered the Korobov space
\begin{equation}
E_a^d = \left\{ f \in L^2([0,1]^d)~:~\exists C > 0~\text{s.t.}~\forall \bm j \in \Z^d,~|c_{\bm j} (f)| \leq C\prod_{k=1}^d (1+|j_k|)^{-a} \right\}
\end{equation} 
where $a > 1$ is a smoothness parameter. In general, $a$ may take on any real value greater than one, but integer values have an interpretation in terms of the order of differentiability~\cite{Papageorgiou_Wozniakowski_2010}. One may directly connect $a \in \{2,3,\dots\}$ back to $H^{\bm n}(\T^d)$ and~\cref{equ:fcoef_bound_multidim} by taking $\bm n = (a-2, \dots, a-2)$. There is also a precedent in the literature for our consideration of anisotropic, rather than isotropic, approximations for functions obeying~\cref{equ:fcoef_bound_multidim}. Authors have recently studied tractability questions in anisotropic Korobov spaces~\cite{Kritzer_Pillichshammer_Wozniakowski_2014,Novak_Wozniakowski_2008,Papageorgiou_Wozniakowski_2010} in addition to approximation in anisotropic Sobolev and Besov spaces~\cite{Griebel_Hamaekers_2014,Sickel_Ullrich_2009}. There is also a long tradition of dimensionally and spatially adaptive refinement within the context of sparse-grid interpolation with a piecewise or Lagrange polynomial basis, e.g.~\cite{Bungartz_Griebel_2004, gunzburger2014stochastic, jakeman2011characterization, jakeman2013minimal, jakeman2012local, khakhutskyy2016spatially, klimke2005algorithm, ma2009adaptive, narayan2014adaptive, Nobile_Tamellini_Tempone_2016, nobile2008anisotropic, pflueger12spatially, pfluger2010spatially, stoyanov2018adaptive, stoyanov2017predicting, Stoyanov_Webster_2016}. Our method differs from previous work by using trigonometric basis functions rather than polynomials and by estimating the anisotropy using a fit to a sharp bound on the decay rate of Fourier coefficients \cref{equ:fcoef_bound_multidim}; see \cref{subsec:anisotropy}. Furthermore, to the authors' knowledge, this method is the first to use adaptive refinement with trigonometric interpolation on sparse grids.

\subsection{Quasi-optimal projection space}
\label{subsec:quasi_opt_projection}

Consider the projection of $f \in H^{\bm n}(\T^d)$ onto a space of real trigonometric polynomials defined by the finite lower complete\footnote{A set $\Lambda$ is called \textit{lower} if $\bm \nu \in \Lambda$ implies $\{ \bm i \in \N^d~:~\bm i \leq \bm \nu \} \subset \Lambda$, where $\bm i \leq \bm \nu$ if and only if $i_k \leq \nu_k$ for each $1 \leq k \leq d$. \label{foot:lower}} multi-index set $\Lambda \subset \N^d$:
\begin{equation} 
\P_{\Lambda} = \text{span} \bigcup_{\bm \nu \in \Lambda} \P_{\bm \nu}
= \text{span} \bigcup_{\bm \nu \in \Lambda} \ \bigotimes_{k=1}^d \P_{\nu_k} \, 
\label{equ:trig_poly_space}
\end{equation} 
where 
\begin{equation} 
\P_n = \text{span}\{ \exp(2 \pi \text{i} \, j \, x) \}_{j=-n}^n \, , \qquad \text{i}^2 = -1 \, , \label{equ:trig_poly_integer}
\end{equation} 
which is the set of all univariate trigonometric polynomials of degree at most $n$. In the context of the sparse grid construction presented in~\cref{sec:sparse_grids}, it is more convenient to represent the spaces in terms of non-negative multi-index sets. As a result, we introduce the re-indexing
\begin{equation}
\phi_\nu (x) = \exp(2 \pi \text{i} \, \sigma(\nu) \, x), \quad \nu \in \N
\end{equation}
where
\begin{equation}
\sigma(\nu) = \begin{cases} 
-\nu/2, & \nu~\text{even} \\
(\nu+1)/2, & \nu~\text{odd}
\end{cases}  \, ,
\label{equ:sigma}
\end{equation}
so that
$$
\P_n = \text{span} \{ \phi_\nu \}_{\nu=0}^{2n} \, .
$$
Note that when using the scalar subscript on $\P_n$ we are in fact referring to the one dimensional space defined by the lower-complete set $\{\nu \in \N : \nu \leq n \}$.

Let $f_{\Lambda}$ be the best approximation to $f$ in $\P_{\Lambda}$ in the $L^2(\T^d)$ sense, i.e.,
$f_{\Lambda}$ is defined from the orthogonal decomposition of $f$ in terms of trigonometric polynomials $\phi_{\bm \nu}(\bm x) \in \P_{\Lambda}$:
\begin{equation}
f_{\Lambda}(\bm x) = \sum_{\bm \nu \in \Lambda} c_{\bm \nu} \, \phi_{\bm \nu}(\bm x), \qquad \phi_{\bm \nu}(\bm x) = \prod_{k=1}^d \phi_{\nu_k}(x_k)
\label{equ:trunc_fourier_series}
\end{equation}
which is the familiar Fourier series. It follows that the best $M$-term approximation space for projection is associated with the $M$ largest Fourier coefficients of $f$. Therefore, taking the upper bound in~\cref{equ:fcoef_bound_multidim} with $\alpha_k = n_k+2$, we obtain the quasi-optimal space
\begin{equation}
\Lambda^{\bm \alpha}(L) = \{ \bm i \in \N^d~:~(\bm i + \bm 1)^{\bm \alpha} \leq L \}, \,
\quad \text{where} \quad \bm \nu^{\bm \alpha} = \prod_{k=1}^d \nu_k^{\alpha_k}.
\label{equ:hyperbolic}
\end{equation}
and the parameter $L \in \N$ discretizes the multi-index space into levels. Note that the structure of the multi-indexes corresponds to a hyperbolic space, which is in contrast to the total-degree space, commonly used for sparse grids~\cite{Barthelmann_Novak_Ritter_2000, Hallatschek_1992, Novak_Ritter_1996, Novak_Ritter_1999}:
\begin{equation}
\Lambda^{\bm \alpha}_{TD}(L) = \{ \bm i \in \N^d~:~\bm \alpha \cdot \bm i \leq L \}
\label{equ:total_degree}
\end{equation}
where $\bm \alpha \cdot \bm i = \sum_{k=1}^d \alpha_k i_k$.

\subsection{Quasi-optimal interpolation}
\label{subsec:quasi_opt_interpolation}

Projection yields the optimal $L^2$ error, but computing $c_{\bm \nu}$ to a sufficient accuracy involves a number of function samples typically much larger than the size of the basis $\P_{\Lambda(L)}$. In contrast, interpolation requires exactly the same number of samples as basis functions at the cost of slight reduction of accuracy. Thus, in many practical situations, interpolation results in a better overall convergence rate.

In this subsection, we consider $f_{\Lambda(L)}$ as an interpolatory (rather than projective) approximation to $f$. Specifically, $f_{\Lambda(L)}$ is obtained from applying the interpolation operator $I_{\Lambda(L)}$ to $f$, where the operator is exact for all functions in $\P_{\Lambda(L)}$. Following classical results in interpolation (e.g.,~\cite{Gautschi_2012}), for all $\phi \in \P_{\Lambda(L)}$ we have
\begin{align*}
\left\| f-f_{\Lambda(L)} \right\|_{\infty} &= \left\| f-\phi + \phi-f_{\Lambda(L)}\right\|_{\infty}  \\
&= \left\| f-\phi + I_{\Lambda(L)}[ \phi-f] \right\|_{\infty} \\
&\leq \left\|f-\phi\right\|_{\infty} +\left\|I_{\Lambda(L)}\right\| \, \left\| f-\phi\right\|_{\infty} \, ,
\end{align*}
where
\begin{equation}
\| I_{\Lambda(L)} \| = \sup_{\|u\|_{\infty}=1} \left\| I_{\Lambda(L)}[u] \right\|_{\infty}
\label{equ:operator_norm}
\end{equation} 
and $\| \cdot \|_{\infty}$ is the $L^\infty$ norm on $\T^d$. Taking the infimum over all $\phi \in \P_{\Lambda(L)}$, we obtain 
\begin{equation}
\| f-f_{\Lambda(L)} \|_{\infty} \leq \left( 1+ \L_{\Lambda(L)} \right) \inf_{\phi \in \P_{\Lambda(L)}} \| f- \phi \|_{\infty} \, ,
\label{equ:interp_error_lebesgue}
\end{equation}
where $\L_{\Lambda(L)}$ is norm of the interpolation operator, commonly called the \textit{Lebesgue constant}. Note that in our context the word \textit{constant} is a misnomer since it strongly depends on the approximation set $\Lambda(L)$ and the specific choice of samples; a detailed discussion is included in \cref{sec:sparse_grids}. Using \cref{equ:trunc_fourier_series}, we observe that
\begin{equation}
\inf_{\phi \in \P_{\Lambda(L)}} \| f-\phi \|_{\infty} \leq \left\| f- \sum_{\nu \in \Lambda(L)} c_{\bm \nu} \, \phi_{\bm \nu} \right\|_{\infty} 
\label{equ:projection_to_Linfty}
\end{equation}
where the coefficients $c_{\bm \nu}$ come from the best $L^2$ approximation. We can chain together \cref{equ:interp_error_lebesgue}-\cref{equ:projection_to_Linfty} and note that the only difference from the optimal $L^2$ approximation comes from using the $L^\infty$ rather than $L^2$ norm. We end up deriving the same quasi-optimal space as in Equation~\cref{equ:hyperbolic} by heuristically approximating interpolation error as a combination of $L^2$ projection error and the Lebesgue constant. Other papers have included more parameters in the quasi-optimal interpolation space in order to incorporate the effects of the Lebesgue constant~\cite{Stoyanov_Webster_2016}, but this is not necessary in our context; see \cref{subsec:lebesgue}.

\subsection{Estimating anisotropy}
\label{subsec:anisotropy}
The specific values of the entries in the anisotropy vector $\bm \alpha$, while critical for constructing a quasi-optimal approximation, are seldom known \textit{a priori}. In this section, we describe a method for estimating the anisotropy from an already constructed approximation $f_{\Lambda(L)}$ for some lower set $\Lambda(L)$. By definition, since $f_{\Lambda(L)} \in \P_{\Lambda(L)}$
$$
f_{\Lambda(L)}(\bm x) = \sum_{\bm \nu \in \Lambda(L)} \hat{c}_{\bm \nu} \, \phi_{\bm \nu}(\bm x)
$$
where $\hat{c}_{\bm \nu}$ are either the projection coefficients from \cref{subsec:quasi_opt_projection} or a corresponding set of interpolation coefficients. Often times, the $\hat{c}_{\bm \nu}$ are explicitly computed as part of the respective projection or interpolation procedure and hence available at no additional cost. If the estimate in~\cref{equ:fcoef_bound_multidim} bounds the decay of $\hat{c}_{\bm \nu}$ sharply, then
\begin{equation} 
|\hat{c}_{\bm \nu}| \approx C(f) \, \prod_{k=1}^d (1+ |\sigma(\nu_k)|)^{-\alpha_k},
 \hspace{2em} \forall \bm \nu \in \Lambda(L)
\label{equ:interp_coef_bound}
\end{equation} 
where we make use of $\sigma(\nu)$ in~\cref{equ:sigma} to re-index $\bm j \in \Z^d$ in the estimate~\cref{equ:fcoef_bound_multidim} to $\bm \nu \in \Lambda(L)$. 

The hyperbolic space defined in ~\cref{equ:hyperbolic} does not depend on the constant $C(f)$; thus, we focus our attention on estimating $\bm \alpha$, and we include the constant only as a regularizing term. If $\Lambda(L)$ is defined by~\cref{equ:hyperbolic} and $L \geq 2$, we can replace the approximate sign in~\cref{equ:interp_coef_bound} by an equal sign and solve the system of equations; however, in practice, the estimate is only an asymptotic upper bound and the individual coefficients $\hat{c}_{\bm \nu}$ can vary in the preasymptotic regime, which gives an effect similar to noise. Thus, we take more samples in each direction and solve for the effective rates of decay from an over-determined set of equations. Taking the log of both sides and changing signs, we obtain
\begin{equation} 
-\log(|\hat{c}_{\bm \nu}|) \approx - \log(C(f)) + \bm \alpha \cdot \log(\bm{\tilde \sigma}(\bm \nu) + \bm 1), \hspace{2em} \forall \bm \nu \in \Lambda(L) .
\label{equ:log_intermediate}
\end{equation} 
where
\begin{equation}
\log(\bm i) = (\log(i_1), \dots, \log(i_d)) \, , \qquad \bm{\tilde \sigma}(\bm \nu) = (|\sigma(\nu_1)|, \dots, |\sigma(\nu_d)|).
\label{equ:sigma_multidim_def}
\end{equation}
We average out the fluctuations in the coefficient values by we taking the least-squares solution, i.e., the solution that minimizes the $\ell^2$ norm
\begin{equation}
\min_{\bm \alpha \in \R^d,\,\bar C \in \R} \ \frac12 \ \sum_{\bm \nu \in \Lambda(L)} (\bar C + \bm \alpha \cdot \log(\bm{\tilde \sigma}(\bm \nu) + \bm 1) + \log(|\hat{c}_{\bm \nu}|))^2 \, ,
\label{equ:anisotropy_lsq}
\end{equation}
which can be written in a matrix form
\begin{equation}
    \min_{\bm v \in \R^{d+1}} \frac{1}{2} \| \bm A \bm v - \bm b \|_2^2,
    \label{equ:anisotropy_lsq2}
\end{equation}
where the rows of $\bm A$ have the form 
$$
\left(1, \, \log(|\sigma(\nu_1)| + 1), \, \log(|\sigma(\nu_2)| + 1), \, \cdots, \, \log(|\sigma(\nu_d)| + 1) \right),
$$
the vector $\bm b$ holds the corresponding entries of $-\log(|\hat c_{\bm \nu}|)$, and the solution vector is $\bm v = (\bar C, \alpha_1, \cdots, \alpha_d)^T$. Note that $\bar C$ is not the same as $\log(C(f))$ from~\cref{equ:fcoef_bound_multidim}, since $C(f)$ defines an upper bound while the curve defined by $\bar C$ and $\bm \alpha$ has coefficients $-\log(\hat c_{\bm \nu})$ both above and below.

The derivation of the anisotropic rates $\bm \alpha$ links them to integers orders of differentiability; however, in our context, we are only interested in the quasi-optimal approximation space and w.l.o.g. we can allow $\bm \alpha$ to be any positive real numbers. We see this by observing that the $L$ parameter is a dummy discretization variable in \cref{equ:hyperbolic} and more obviously in \cref{equ:total_degree}. If we multiply both $\bm \alpha$ and $L$ by the same positive constant, the total degree space \cref{equ:total_degree} remains the same; similarly, in the hyperbolic case \cref{equ:hyperbolic}, we can multiply $\bm \alpha$ and raise $L$ to the same positive power. For the isotropic total degree space with $\bm \alpha = \bm 1$, $L$ indicates the total degree polynomial order; however, the relation is lost for an isotropic space with $\bm \alpha \neq \bm 1$ and even more so for an anisotropic space. In the hyperbolic cross section construction, the notion of $L$ does not relate to polynomial order even for $\bm \alpha = \bm 1$. In both cases, the factor that determines the structure of the quasi-optimal space is the ratio between the pairs of coefficients in $\bm \alpha$. Since most least-squares solvers operate on real valued numbers, it is convenient to remove any integer restrictions in \cref{equ:anisotropy_lsq}-\cref{equ:anisotropy_lsq2} and seek $\bm v \in \mathbb{R}^{d + 1}$.

Equation \cref{equ:anisotropy_lsq2} admits a unique solution so long as $\bm A$ has full column rank, i.e., so long as we have at least two coefficients in each direction to estimate the corresponding decay rate. However, since the coefficients may not decay monotonically and since the accuracy of the solution heavily depends on the condition number of the matrix $\bm A$, the approximation with only two coefficients will not suffice. The estimated $\bm \alpha$ may be too inaccurate or even yield negative decay rates, which according to \cref{equ:hyperbolic} results in $\Lambda(L)$ with infinitely many multi-indexes. Nevertheless, we employ the estimate in an adaptive refinement strategy presented in \cref{subsec:adaptive_refinement}, and in \cref{rem:adhoc_stabilization} we propose an ad-hoc strategy, specific to the refinement procedure, that would allow us to move forward with the adaptive steps even if some of the computed $\alpha_k$ are negative.

\section{Sparse trigonometric interpolation}
\label{sec:sparse_grids}

In this section, we describe a sparse-grid approach for constructing an interpolant $f_{\Lambda(L)} \in \P_{\Lambda(L)}$ based on the values of $f$ at a set of nodes $\bm x_1, \dots, \bm x_M \in \T^d$. Here, we make no optimality assumptions about $\Lambda(L)$; we only assume that $\Lambda(L)$ is a lower set. After discussing the sparse interpolation algorithm, we examine the sparse-grid Lebesgue constant. We conclude by coupling the anisotropy-estimation procedure of \cref{subsec:anisotropy} with a general sparse trigonometric interpolation algorithm in order to produce a dimensionally adaptive interpolant.

\subsection{One-dimensional and fully tensorized rule}
\label{subsec:1d_fulltensor_grid}
First, we define the one dimensional $m$-point trigonometric interpolation rule. The nodes and (global) basis functions are
\begin{equation} 
x_j = \frac{j}{m}, \hspace{2em} \phi_j(x) = \exp(2\pi \text{i} \, \sigma(j) \, x), \hspace{2em} j=0,1,\dots,m-1,
\label{equ:rule_1d}
\end{equation}
Because approximations using cosines and sines up to mode $n$ require complex exponentials with powers ranging between $-n \leq j \leq n$~\cite{Stoer_Bulirsch_1993}, we use the re-indexing $\sigma(j)$ first introduced in~\cref{equ:sigma}. To resolve all sine and cosine modes up to $n$, it is necessary and sufficient to take $m = 2n+1$ points~\cite{Stoer_Bulirsch_1993}. The interpolation operator is
\begin{equation} 
\mathcal{U}^{m} : C^0(\T) \to \P_{ (m-1)/2 }, \hspace{2em} \mathcal{U}^{m}[f](x) = \sum_{j=0}^{m-1} \hat{c}_j \, \phi_j(x),
\label{equ:interpolant_1d}
\end{equation} 
where $\hat{c}_j$ are the interpolation coefficients, $\P_n$ is defined by~\cref{equ:trig_poly_integer}, and $C^0(\T)$ is the set of all real-valued continuous functions on $\T$.

\begin{remark}
\label{rem:odd_points}
One-dimensional (non-adaptive) trigonometric interpolation with an even number of points has precedent in the literature~\cite{Ehlich_Zeller_1966, Stoer_Bulirsch_1993}. Indeed, much work with Fourier transforms on sparse grids uses $2^l$ points at each level since powers of 2 are highly amenable to fast Fourier transforms~\cite{Griebel_Hamaekers_2014,Hallatschek_1992}. However, for the purposes of interpolation, an even number of points $m=2n$ exactly reproduces all modes up to $n-1$ due to the missing conjugate exponent for mode $n$. Put differently, there is an additional basis function for $m=2n$ vs. $m=2n-1$ without a general increase in exactness in terms of both sines and cosines. Therefore, we use an odd number of interpolation nodes, which means $\mathcal{U}^{m}$ is exact up to mode $(m-1)/2$. 
\end{remark} 

The interpolation conditions at the nodes $x_j$
$$
\mathcal{U}^{m}[f](x_j) = f(x_j), \hspace{2em} j=0,1,\dots,m-1,
$$
can be imposed by selecting the coefficients as~\cite{Griebel_Hamaekers_2014,Hallatschek_1992, Stoer_Bulirsch_1993}
\begin{equation} 
\hat{c}_j = \frac{1}{m} \sum_{p=0}^{m-1} f(x_p) \left[ \phi_j(x_p) \right]^*
\label{equ:dft_1d}
\end{equation} 
where $[\, \cdot \, ]^*$ denotes complex conjugation. With some algebra work, it becomes apparent that~\cref{equ:dft_1d} is a normalized and re-indexed discrete Fourier transform. Because the target function $f \in C^0(\T)$ is real-valued, symmetry in $\hat{c}_j$ and $\phi_j$ makes the interpolant real-valued in exact arithmetic; in implementation when computing~\cref{equ:interpolant_1d} numerically, we only compute the real part.

We extend the one dimensional construction to a multidimensional context using tensor products of the points and basis functions expressed in multi-index notation. Let $\bm m = (m_1, m_2, \dots, m_d)$ represent the vector with (potentially) different number of points in each dimension, then
\begin{align*}
    \bm x_{\bm j} &= \left( x_{j_1}, x_{j_2}, \cdots, x_{j_d} \right ), \\
    \phi_{\bm j}(\bm x) &= \prod_{k=1}^d \phi_{j_k}(x_k) 
                = \exp\left(2\pi \text{i} \, \sum_{k=1}^d \sigma(j_k) \cdot x_k\right).
\end{align*}
The anisotropic fully tensorized operator becomes:
\begin{align*}
\mathcal{U}^{\bm m} : C^0(\T^d) & \to \bigotimes_{k=1}^d \P_{(m_k-1)/2}\, ,  \\
\mathcal{U}^{\bm m}[f](\bm x) &= \sum_{\bm j \leq \bm m-\bm 1} \hat{c}_{\bm j} \phi_{\bm j}(\bm x) \, .
\end{align*} 
As one may expect, the fully tensorized interpolation coefficients are analogous to the one-dimensional case:
\begin{equation} 
\hat{c}_{\bm j} = \frac{1}{m_1 \cdots m_d} \ \ \sum_{\bm p \leq \bm m-\bm 1} f(\bm x_{\bm p}) \left[ \phi_{\bm j}(\bm x_{\bm p})\right]^* \, ,
\label{equ:dft_fulltensor}
\end{equation} 
which is a normalized and re-indexed $d$-dimensional discrete Fourier transform.

The fully tensorized construction is easy to implement, using only one-dimensional nodes and basis functions and employing a suitable algorithm for fast-Fourier transform. However, the resulting approximation belongs to a fully tensorized space which is very far from optimal. First, we observe that the extra point incurred by the even rules in one-dimension combines with all the ``good'' points in other dimensions and results in a much higher penalty, e.g., in $6$-dimensions the $4$-point rule has $4096$ points and it covers the same basis as the $3$ point rule with only $729$ points, thus wasting the majority of the computational effort. Restricting our attention to rules with only odd number of points, we then look at our estimate for the quasi-optimal basis. At $L = 2$, in $6$ dimensions, the space has $7$ exponential powers which results in $13$ basis functions (all powers except zero require two basis functions), the smallest fully tensorized space including $L=2$ hyperbolic space has the aforementioned $729$ points. Fully tensorized interpolation is not a feasible approach in a multidimensional context.

\subsection{Sparse-grid interpolation}
\label{subsec:general_sparse_grid}

Sparse-grid interpolation aims at exploiting the implementational simplicity of fully tensorized rules while alleviating (and sometimes completely avoiding) the restrictions on the basis space. To this end, sparse-grid algorithms employ a family of one dimensional interpolation rules with different number of points and basis functions and combine (superimpose) a set of anisotropic full tensors interpolants into a single grid. The set of tensors is chosen so that the combined approximation space includes a desired (quasi-) optimal space with as little extra basis functions as possible.

Starting with the one-dimensional nodes and basis~\cref{equ:rule_1d}, we select the family of rules via the node growth $m(l)$; i.e., $m(l)$ is a strictly increasing function indicating the number of points on level $l \geq 0$. See~\cref{rem:nested_points} for the specific choice used in our examples. Following the approach used in~\cite{Stoyanov_Webster_2016}, we define the surplus operators
$$
\Delta^{m(l)} = \mathcal{U}^{m(l)}-\mathcal{U}^{m(l-1)}, \hspace{2em} \Delta^{\bm m(\bm i)} = \bigotimes_{k=1}^d \Delta^{m(i_k)},
$$
with the convention $\Delta^{m(0)} = \mathcal{U}^{m(0)}$. For any lower set $\Theta(L)$, we define the generalized interpolation operator
\begin{equation}
I_{\Theta(L)} = \sum_{\bm i \in \Theta(L)} \Delta^{\bm m(\bm i)}.
\label{equ:sparse_operator_delta}
\end{equation}
Our objective is to relate $\Theta(L)$ to the quasi-optimal $\Lambda(L)$, but first we observe that
\begin{equation} 
\mathcal{U}^{\bm m(\bm i)} = \sum_{\bm j \leq \bm i} \Delta^{\bm m(\bm j)}
\label{equ:fulltensor_decomposed}
\end{equation} 
since the sum is telescoping~\cite{Bungartz_Griebel_2004, nobile2008sparse, Stoyanov_Webster_2016}. We wish to express 
\begin{equation}
I_{\Theta(L)} = \sum_{\bm i \in \Theta(L)} t_{\bm i} \, \mathcal{U}^{\bm m(\bm i)}
\label{equ:sparse_fulltensor}
\end{equation}
for some coefficients $\{ t_{\bm i} \}_{\bm i \in \Theta(L)}$. By substituting~\cref{equ:fulltensor_decomposed} into~\cref{equ:sparse_fulltensor} and equating the coefficients of the $\Delta^{\bm m(\bm i)}$ operators, we derive the system
\begin{equation} 
\sum_{\substack{\bm i \in \Theta(L) \\ \bm j \leq \bm i }} t_{\bm i} = 1, \hspace{2em} \forall \bm j \in \Theta(L) \, .
\label{equ:tj_system}
\end{equation} 
The system~\cref{equ:tj_system} can be expressed as an upper triangular matrix of zeroes and ones, with a diagonal of all ones, so a unique integer solution $\{ t_{\bm i }\}_{\bm i \in \Theta(L)}$ does indeed exist.

The interpolation nodes associated with $I_{\Theta(L)}$ are the union of the nodes of all tensors $ \mathcal{U}^{\bm m(\bm i)} $.
If we want to minimize the number of nodes (and the associated expensive simulation of the target model), it is best to reuse the nodes as much as possible; i.e., we want the nodes associated with $ \mathcal{U}^{m(l)} $ to be a subset of the nodes of $ \mathcal{U}^{m(l+1)} $. It is well known that such nested rules are advantageous for sparse-grid methods, and thus, we restrict out attention to only those $m(l)$ that satisfy the nested property; see~\cref{rem:nested_points}.

Let $\Theta_m(L)$ denote the multi-indexes of the interpolation nodes for the interpolant defined by $\Theta(L)$. Then exploiting the nested structure of the rule gives
\begin{equation}
\Theta_m(L) = \bigcup_{\bm i \in \Theta(L)} \left\{ \bm j \in \N^d~:~\bm j \leq \bm m(\bm i)-\bm 1 \right\} \, ,
\label{equ:theta_m}
\end{equation}
which comes from
$$
\{ \bm x_{\bm j} \}_{\bm j \in \Theta_m(L)} = \bigcup_{\bm i \in \Theta(L)} \{ \bm x_{\bm j} \}_{\bm j \leq \bm m(\bm i)-\bm 1} \, .
$$
With~\cref{equ:sparse_fulltensor}-\cref{equ:theta_m}, we can explicitly write the sparse trigonometric interpolant as
\begin{equation}
I_{\Theta(L)}[f](\bm x) = \sum_{\bm j \in \Theta_m(L)} \ \sum_{\substack{\bm i \in \Theta(L) \\ \bm j \leq \bm m(\bm i)-\bm 1}} t_{\bm i} \, \hat{c}^{\bm i}_{\bm j} \, \phi_{\bm j}(\bm x) = \sum_{\bm j \in \Theta_m(L)} w_{\bm j} \, \phi_{\bm j}(\bm x) \, .
\label{equ:sparse_interpolant}
\end{equation}
Here, $\hat{c}^{\bm i}_{\bm j}$ and $\phi_{\bm j}$ are defined in~\cref{subsec:1d_fulltensor_grid}, and
\begin{equation}
w_{\bm j} = \sum_{\substack{\bm i \in \Theta(L) \\ \bm j \leq \bm m(\bm i)-\bm 1}} t_{\bm i} \, \hat{c}^{\bm i}_{\bm j}\, , \hspace{3em} \bm j \in \Theta_m(L).
\label{equ:wj}
\end{equation}
Since each of the tensor operators exactly reproduces the basis functions, constructing $I_{\Theta(L)}[\phi_{\bm j}](\bm x)$ means that the only non-zero coefficient will be the corresponding $ \hat c_{\bm j}^{\bm i} $ and from~\cref{equ:tj_system} follows that the interpolant is exact for all basis functions, i.e., the union of the space of all tensors. Therefore, Theorem 1 in~\cite{Stoyanov_Webster_2016} applies; since the exactness of level $l \in \N$ is $(m(l)-1)/2$, then for an arbitrary lower trigonometric polynomial space $\P_{\Lambda(L)}$, we define the optimal sparse  grid tensor set by
\begin{equation}
\Theta^{opt}(L) = \{ \bm i \in \N^d~:~(\bm m(\bm i-\bm 1)+\bm 1)/2 \in \Lambda(L) \} \, .
\label{equ:theta_opt}
\end{equation}
That is, $\Theta^{opt}(L)$ is the smallest set of tensors that results in an interpolant which is exact for $\P_{\Lambda(L)}$. Depending on the choice of $m(l)$, the actual interpolation space may be larger.

\begin{remark}
\label{rem:nested_points} 
We choose $m(l)=3^l$, which gives us nodes that are both nested and odd, i.e., we avoid the cost of extra basis functions noted in~\cref{rem:odd_points}. The Radix-3 FFT algorithms are slightly less efficient than the Radix-2 and Radix-4 variants; however, from~\cref{equ:wj}, we observe that the weight can be pre-computed and reused every time we need to compute the value of the interpolant. Thus, the FFT procedure is a one-time effort resulting in a small increase in computational cost, which is far offset by the reduction in interpolation nodes and model simulations. Furthermore, when targeting hyperbolic cross-section interpolation space, e.g.~\cref{equ:hyperbolic}, exponentially growing $m(l)$ have a natural advantage. Suppose that $\Theta(L)$ is chosen as a total degree multi-index space:
\begin{equation}
\Theta(L) = \left\{ \bm i \in \N^d : \sum_{k=1}^d \alpha_k i_k \leq L \right\} \quad \Rightarrow \quad \prod_{k=1}^d \left( 3^{i_k} \right)^{\alpha_k} \leq 3^L \, .
\label{equ:theta_total}
\end{equation}
Then from the definition of $\Theta_m(L)$ in~\cref{equ:theta_m} and the one-to-one correspondence between nodes and basis functions, we have that the resulting space will include all $\phi_{\bm j}$ for $\bm j$ such that
$$
\bm j + \bm 1 \leq \bm m(\bm i), \, \quad \text{i.e.} \quad 
j_k + 1 \leq 3^{i_k} \quad \forall 1 \leq k \leq d.
$$
Combining the above with~\cref{equ:theta_total}, we get the anisotropic cross-section space:
$$
\prod_{k=1}^d (j_k + 1)^{\alpha_k} \leq 3^L \, .
$$
Therefore, using an exponential $m(l)$ and total degree multi-index selection produces an interpolant in a hyperbolic cross-section space, modulo a constant in the exponent and some rounding in the actual implementation due to non-integer anisotropic rates.

\end{remark}

\subsection{Lebesgue constant}
\label{subsec:lebesgue}

Consider the Lebesgue constant $\L_n$ for $(2n+1)$-point trigonometric interpolation. Well-known estimates for $\L_n$ exist~\cite{Ehlich_Zeller_1966}, including quite sharp ones~\cite{Dzyadyk_Dzyadyk_Prypik_1981,Rivlin_1974}. In terms of the Lebesgue constant, trigonometric interpolation is closely related to polynomial interpolation with the Clenshaw--Curtis and Chebyshev nodes~\cite{Clenshaw_Curtis_1960, Ehlich_Zeller_1966}. In~\cite{Rivlin_1974}, Rivlin showed the equality
\begin{equation}
\L_n = \frac{2}{\pi} \ln(n) + \beta_n, \hspace{2em} n \geq 1,
\label{equ:lebesgue_trig_exact}
\end{equation}
where $\beta_n$ decreases monotonically from $5/3$ to
$$
\frac{2}{\pi} \left( \ln\left( \frac{16}{\pi} \right) + \gamma \right) \approx 1.404,
$$
where $\gamma$ is the Euler--Mascheroni constant. For our family of one dimensional rules defined by $m(l) = 3^l$, the Lebesgue constant grows as a logarithm in the number of points and linearly in level.

In general, no sharp estimates exist for the Lebesgue constant of sparse interpolation with space of exactness $\P_{\Lambda(L)}$. However, Lemma 3.1 in~\cite{Chkifa_Cohen_Schwab_2014} yields an upper bound, i.e., when the Lebesgue constant of the one-dimensional family of rules exhibits a polynomial growth, the sparse grids constant grows no faster than the polynomial with one additional power. Let $\Theta(L)$ be the optimal multi-index set corresponding to $\Lambda(L)$. Because the one-dimensional rules obey~\cref{equ:lebesgue_trig_exact}, we can bound the Lebesgue constant of the sparse interpolation operator $I_{\Theta(L)}$ in~\cref{equ:sparse_interpolant} as
\begin{equation}
\left\| I_{\Theta(L)} \right\| \leq C^d \, (\#\Theta(L))^2 \, ,
\label{equ:lebesgue_sparse}
\end{equation}
where $\#\Theta(L)$ is the number of multi-indexes in $\Theta(L)$.
As in~\cite{Chkifa_Cohen_Schwab_2014}, the operator norm is defined by~\cref{equ:operator_norm}.

In~\cite{Stoyanov_Webster_2016}, Stoyanov and Webster considered many different one-dimensional rules, some of which have Lebesgue constants that grow slowly at first but increase quite rapidly after only a few levels. As a result, the authors included an additional parameter to account for the different ``effective'' Lebesgue constant in each dimension since some dimensions may have far fewer points than others. However, we do not consider such a correction in our context for the following reasons. First, the rapid increase in Lebesgue constants in~\cite{Stoyanov_Webster_2016} was not observed for the Clenshaw--Curtis nodes, which have the same linear-in-level and logarithmic-in-nodes growth as our rule. Second, the variability of $\beta_n$ in~\cref{equ:lebesgue_trig_exact} is small. Third, by using~\cref{equ:lebesgue_trig_exact} in~\cref{equ:interp_error_lebesgue} and following the derivation in~\cref{subsec:anisotropy}, the Lebesgue penalty term in~\cref{equ:anisotropy_lsq} is the logarithm of a logarithm, which is negligible in practical situations.

\subsection{Adaptive refinement}
\label{subsec:adaptive_refinement}

In this subsection, we will link the anisotropy least-squares problem from~\cref{subsec:anisotropy} to the sparse trigonometric interpolation algorithm in~\cref{subsec:general_sparse_grid}. We use~\cref{equ:anisotropy_lsq}, but with the sparse interpolation coefficients $w_{\bm j}$ given by~\cref{equ:wj}, which are linear combinations of $\hat{c}^{\bm i}_{\bm j}$ over the constituent tensors of the sparse grid. Also, by~\cref{equ:rule_1d,equ:dft_1d}, the sparse discrete Fourier coefficient $w_{\bm j}$ corresponds to mode $\bm \sigma(\bm j)$. In general, multi-index space is not easily discretizable in the sense of~\cref{equ:hyperbolic} after adaptive refinement takes place, so we omit the dependence on $L$ here. Thus, for a general lower set $\Lambda$ and the corresponding $\Theta_{m}$, the relevant least-squares problem becomes
\begin{equation}
\min_{\bm \alpha \in \R^d,\,\bar C \in \R} \ \frac12 \, \sum_{\bm j \in \Theta_m} (\bar C + \bm \alpha \cdot \log(\bm{\tilde \sigma}(\bm j) + \bm 1) + \log(|w_{\bm j}|))^2 \, .
\label{equ:sparse_lsq}
\end{equation}
where $\Theta_m$ is defined by~\cref{equ:theta_m} and $\bm{\tilde \sigma}(\bm j)$ by~\cref{equ:sigma_multidim_def}.

\begin{figure}
    \centering
    \includegraphics[width=\textwidth]{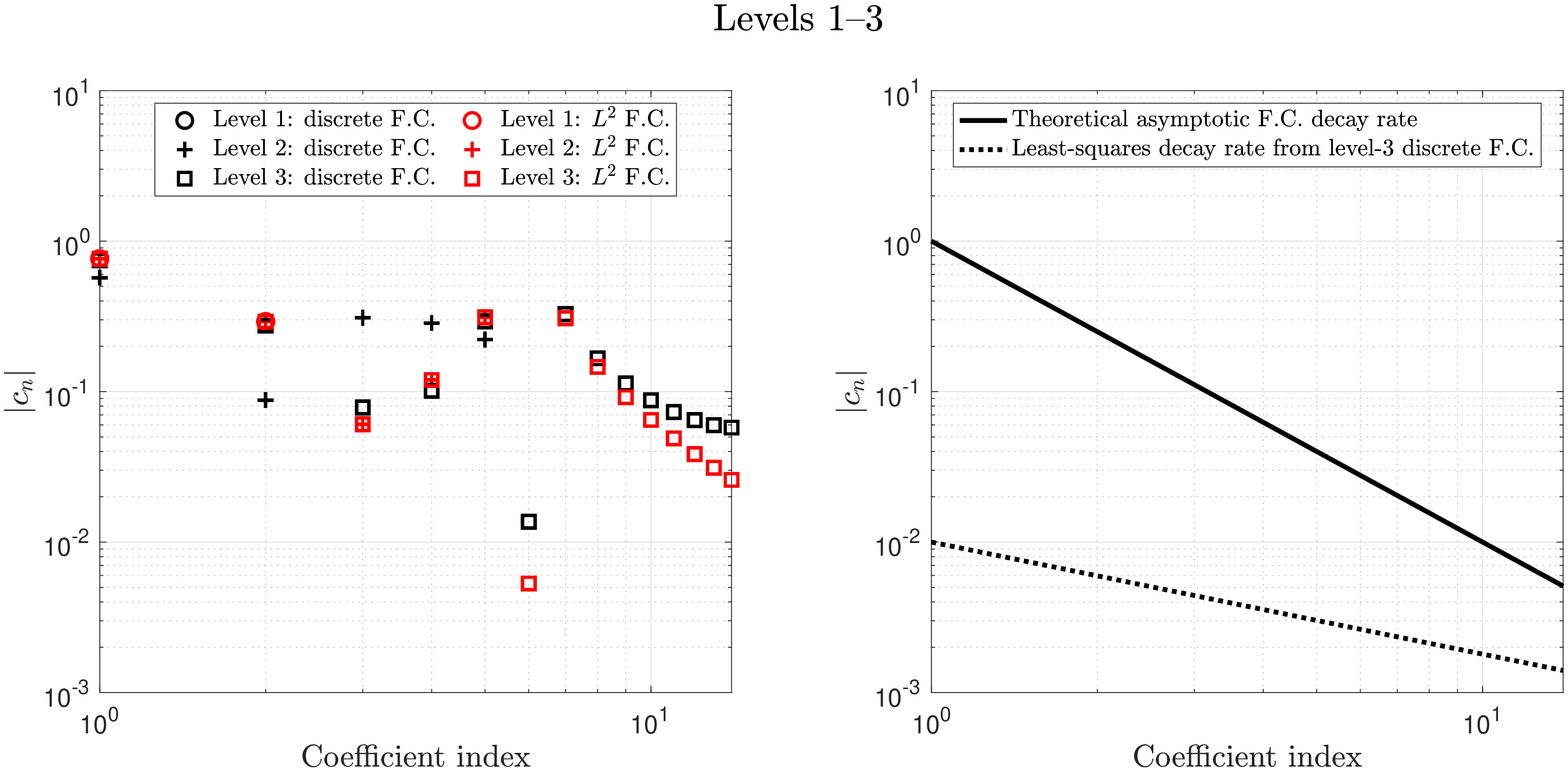}
    \includegraphics[width=\textwidth]{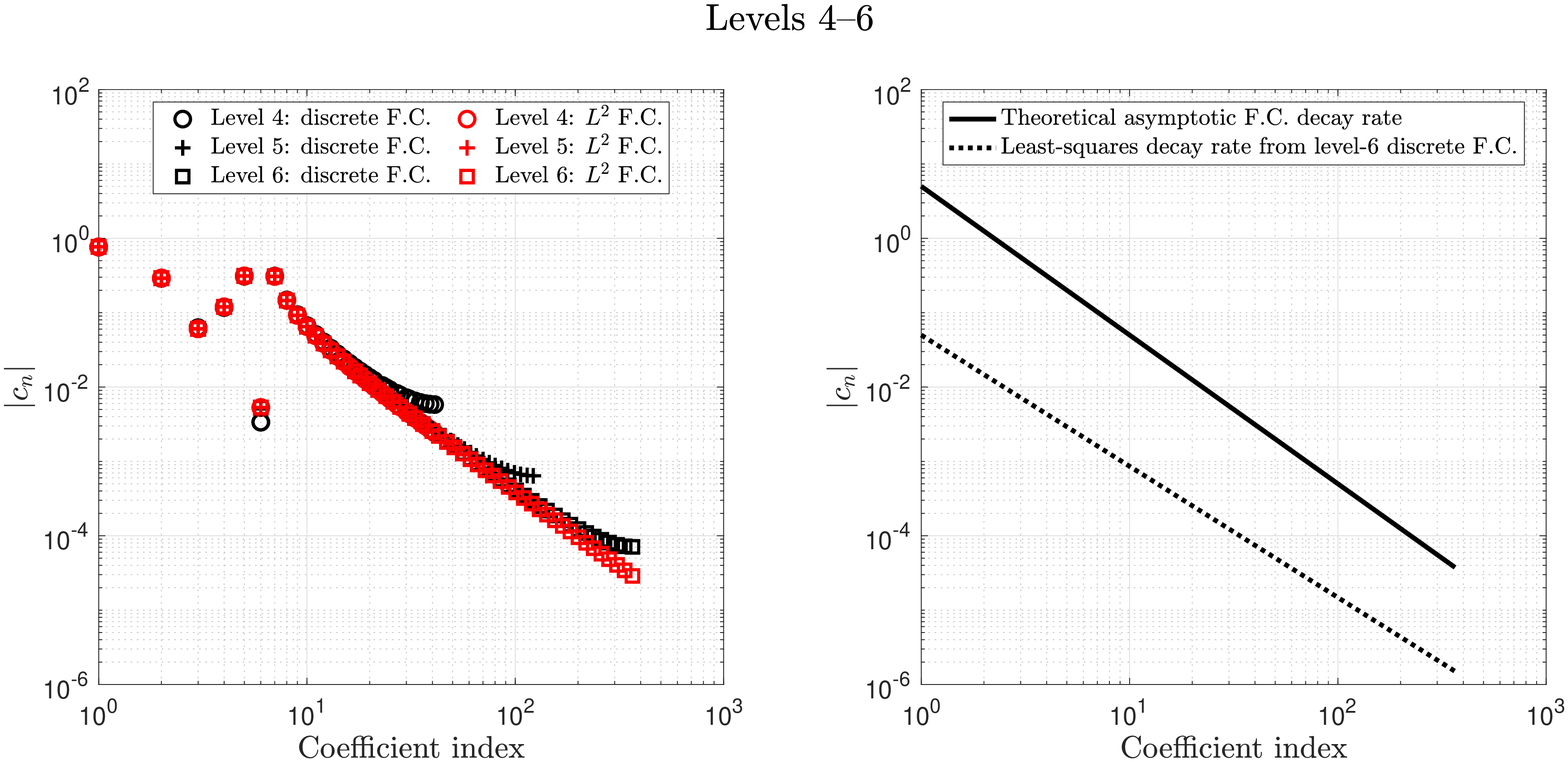}
    \caption{Left column: discrete and $L^2$ Fourier coefficients (F.C.) of \cref{equ:1d_lsq_example_fun} on 1D grids of various sizes. Right column: decay rates from \cref{equ:fcoef_bound_1d} and \cref{equ:sparse_lsq}.}
    \label{fig:fcoef_plot_levels}
\end{figure}

The motivation for the least-squares fitting~\cref{equ:anisotropy_lsq} can be demonstated in the following one-dimensional example
\begin{equation}
f(x) = x \sin(\pi x) + x \sin(5 \pi x), \qquad x \in [-1,1] ,
\label{equ:1d_lsq_example_fun}
\end{equation}
where $f \in H^0(\T)$. In~\cref{fig:fcoef_plot_levels}, we show the computed continuous and discrete Fourier coefficients, the theoretical decay rate according to~\cref{equ:fcoef_bound_multidim}, and the decay rate estimate coming from~\cref{equ:anisotropy_lsq}. We observe that on each level, for the largest indexes in $\Lambda$, the discrete Fourier coefficients ($\circ$) differ systematically from the continuous Fourier coefficients ($+$). This can be explained by observing that the discrete Fourier transform~\cref{equ:dft_1d} is a left-hand Riemann sum discretization of~\cref{equ:fcoef_def}; therefore, the largest indexes in $\Lambda$ correspond to the highest frequencies, and the discretization is not able to resolve those to the same degree of accuracy. A refinement criterion could be based on the surplus or the correction introduced by the high frequencies, e.g., similar to the greedy knapsack problem~\cite{Bungartz_Griebel_2004}, but such refinement would be guided by the least accurate coefficients. This phenomenon is not present in the methods using hierarchical Lagrange approximation where adding more indexes to $\Lambda$ would not alter the current set of polynomial coefficients. Furthermore, there are unpredictable fluctuations in the preasymptotic low-frequencies. Since the breaking point between the two regimes is unknown, we use the least-squares approach defined in~\cref{equ:sparse_lsq} to incorporate {\em all} coefficients and balance out the ``noise-like'' effects.

We show pseudocode for our algorithm in~\cref{alg:adaptive_loop}. Importantly, since the solution of the least-squares problem~\cref{equ:sparse_lsq} is heavily dependent on $\Theta_m$, then one should choose $\Lambda_0$ so that it contains enough points to compute an initial anisotropy estimate that is reliable for each direction. Thus, we select isotropic $\Lambda_0$ and in practical application, the choice is usually guided by considerations regarding the minimum number of samples needed to saturate the computational resources, e.g., the number of computing nodes, see \cref{rem:min_growth} and \cref{rem:adhoc_stabilization}. After we compute an approximate ansitropic coefficients, we want the set $\Lambda^{\bm \alpha}(L)$ in~\cref{equ:hyperbolic} to be invariant when $\bm \alpha$ is multiplied by a constant factor; for instance, $\Lambda^{(1,1)}(L)$ should be equal to $\Lambda^{(2,2)}(L)$, and $\Lambda^{(1,2)}(L)$ equal to $\Lambda^{(2,4)}(L)$. Additionally, we want $\Lambda(L)$ to grow conservatively, not suddenly or dramatically, when incrementing $L$. Thus, we normalize $\bm \alpha$ by dividing by the smallest element before constructing $\Lambda^{\bm \alpha}(L)$.

\begin{algorithm}
\caption{Adaptive refinement algorithm} 
\label{alg:adaptive_loop}
    \hspace*{\algorithmicindent} \textbf{Input:} $L_0 \geq 2$, $f \in H^{\bm n}(\T^d)$, $m(l)$ defined in~\cref{rem:odd_points} \\
    \hspace*{\algorithmicindent} \textbf{Output:} quasi-optimal $\Lambda$ 
\begin{algorithmic}[1]
	\STATE $n \gets 0$
	\STATE Start with isotropic $\Lambda = \Lambda^{\bm 1}(L_0)$ from~\cref{equ:hyperbolic}; define $\Theta$ according to~\cref{equ:theta_opt} \label{line:theta_def_1}
	\STATE Compute the samples of $f$ and load the values into the grid
	\WHILE{$\texttt{num\_samples} < \texttt{budget}$} \label{line:adaptive_termination}
		\STATE Solve~\cref{equ:sparse_lsq} for $\bm{\hat \alpha}$; set $\bm{\hat \alpha} \gets \bm {\hat \alpha} / (\min_k \hat{\alpha}_k$) and apply~\cref{rem:adhoc_stabilization}
		\STATE Find $L_{n+1}$ such that $\Lambda^{\bm{\hat \alpha}}(L_{n+1}) \not\subseteq \Lambda$; define $\Theta^{\bm{\hat \alpha}}(L_{n+1})$ by~\cref{equ:theta_opt} \label{line:theta_def_2} 
		\STATE $\Lambda \gets \Lambda \cup \Lambda^{\bm{\hat \alpha}}(L_{n+1})$;   $\Theta \gets \Theta \cup \Theta^{\bm{\hat\alpha}}(L_{n+1})$;   $n \gets n+1$
		\STATE Compute the samples of $f$ at the new points
	\ENDWHILE 
\end{algorithmic} 
\end{algorithm} 

\begin{remark}
\label{rem:min_growth}
One may modify~\cref{line:theta_def_2} of~\cref{alg:adaptive_loop} so that the number of new points is large enough to exploit parallel computations of the sample of $f$. When the model is sufficiently complex, parallelism must be exploited to make the problem feasible. A possible drawback is that, if the initial grid is too coarse, \cref{alg:adaptive_loop} may add a large number of unnecessary nodes due to an unreliable initial anisotropy estimate.
\end{remark}

\begin{remark}
\label{rem:adaptive_termination}
The quasi-optimal approach aims to construct the best $M$-term approximation for some $M$ given beforehand. For this reason, the termination criterion in~\cref{line:adaptive_termination} relies on reaching some predetermined computational budget. Alternatively, if computing resources are sufficiently abundant, one could terminate upon reaching a desired error tolerance. 
\end{remark}

\begin{remark}
\label{rem:adhoc_stabilization}
Given a black-box model, it is not feasible to determine {\em a priori} the appropriate size of $\Lambda(L_0)$ that would yield a stable initial estimate of the anisotropic coefficients. However, the weights are only used to guide the refinement process. Thus, if we encounter a negative weight $\alpha_k \leq 0$ for some direction $k$, we can simply replace that weight with the smallest positive one, which will force the refinement to put additional points in direction $k$, which in turn will improve the estimate in the following iterations. If all weights are negative, then we continue the refinement using isotropic weights $\bm \alpha = \bm 1$. The correction strategy will allow us to work past negative weights, but it is still possible for a coarse grid to yield positive yet incorrect weight that would deteriorate the convergence. However, the theoretical estimates are only asymptotic and in our numerical examples we observe the opposite behavior, namely that the preasymptotic weights improve the initial error compared to the optimal analytic weights, e.g. in \cref{fig:polynomial_results}. Therefore, in our examples we use isotropic initial $\Lambda(L_0)$ with $L_0 = 3$ which is one more than the absolute minimum.
\end{remark}

\section{Numerical results}
\label{sec:numerical_results}

We include several examples in this section to illustrate the performance of~\cref{alg:adaptive_loop}. We will apply our algorithm to purpose-built periodic polynomials of known anisotropy and then to the chemistry problem that motivated this work. These simulations use the open-source Tasmanian package developed at Oak Ridge National Laboratory~\cite{Tasmanian}, which implements~\cref{alg:adaptive_loop} for sparse trigonometric interpolation.

First, to obtain a theoretical convergence rate for our interpolation algorithm, let $f \in H^{\bm n}(\T)$. Using a theorem of Jackson~\cite{Pinkus_2003}, we can bound the infimum term in \cref{equ:interp_error_lebesgue} by
$$
\inf_{T \in \P_{\Lambda}} \| f - T\|_{\infty} \leq \frac{C(f)}{N^{M+1}}
$$
where $C>0$ is a constant depending on $f$, $M=\min_{k} n_k$, and $N = \# \Theta^{opt}_m$ is the number of nodes. Then, using \cref{equ:lebesgue_sparse} and heuristically approximating $\# \Theta^{opt}$ as $\log(N)$ in light of \cref{rem:nested_points},
\begin{equation}
\| f - I_{\Theta}[f] \|_{\infty} \leq O \left( \log^2(N)/ N^{M+1} \right) 
\label{equ:sparse_convergence_rate}
\end{equation}
for $N$ sufficiently large.

\begin{remark}
\label{rem:total_degree}
(Alternative Function Space) As noted in~\cref{subsec:quasi_opt_interpolation}, much early work on sparse grids sought to approximate function spaces of some total degree~\cref{equ:total_degree}. In Fourier interpolation, the total-degree space is suitable for target functions $f$ having a holomorphic extension in component $k$ within a polyellipse of radius $\alpha_k$ around the real axis, for each $1 \leq k \leq d$. To see why, suppose $f$ is a function satisfying the previous analyticity assumptions. From, e.g.~\cite[p.~27]{Katznelson_2004}, the Fourier coefficients obey the asymptotically sharp estimate
\begin{equation}
|c_{\bm j}(f)| \leq C(f) \, \exp(-\bm \alpha \cdot \bm j) \, , \qquad \bm j \in \Z^d \, .
\label{equ:fcoef_analytic}
\end{equation}
By taking the negative logarithm of the right-hand side of \cref{equ:fcoef_analytic} and ignoring the constant, we obtain the total-degree space \cref{equ:total_degree}. For functions of this type, the least-squares problem~\cref{equ:sparse_lsq} in~\cref{subsec:adaptive_refinement} becomes
\begin{equation}
\min_{\bm \alpha \in \R^d,\,\bar C \in \R} \ \frac12 \, \sum_{\bm j \in \Theta_m} (\bar C + \bm \alpha \cdot \bm{\tilde \sigma}(\bm j) + \log(|w_{\bm j}|))^2 \, ,
\label{equ:sparse_lsq_td}
\end{equation}
where $\bm{\tilde \sigma}(\bm j)$ is defined by~\cref{equ:sigma_multidim_def}. Our numerical examples will include a modification of~\cref{alg:adaptive_loop} that uses the total-degree space~\cref{equ:total_degree} and the least-squares problem~\cref{equ:sparse_lsq_td}.
\end{remark}

\subsection{Periodic polynomials}
\label{subsec:results_polynomials}

We manufacture some multidimensional target functions that are engineered to have a certain order of differentiability and periodicity. We define the univariate functions $g_i : [-1,1] \to \R$ as
\begin{align*}
g_1(x) &= x^3-x \, , \\
g_2(x) &= \frac{x^4}{4} - \frac{x^2}{2} \, ,\\
g_3(x) &= \frac{x^5}{20} - \frac{x^3}{6} + \frac{7x}{60} \, , \\
g_4(x) &= \frac{x^6}{120} - \frac{x^4}{24}  + \frac{7x^2}{120} \, , \\
g_5(x) &= \frac{x^7}{840} - \frac{x^5}{120} + \frac{7x^3}{360}- \frac{31x}{2520} \, ,
\end{align*}
which we have derived by starting with $g_1(x)$ and integrating repeatedly and choosing the constant to preserve periodicity. By construction, $g_k \in H^k([-1,1])$, where we translate $[-1, 1]$ to $[0, 1]$ using a linear transformation and note that the $k+1$-th derivative is discontinuous across the periodic boundary. For $1 \leq i \leq 5$, we normalize in the sup-norm by taking $h_i = g_i / {\|g_i\|_{L^\infty([-1,1])}}$. Thus, the multivariate target functions are
\begin{equation}
f_{\bm i}(\bm x) = \prod_{k=1}^d h_{i_k}(x_k), \hspace{2em} \bm 1 \leq \bm i \leq \bm 5 \, .
\label{equ:multidim_periodic_poly} 
\end{equation}

\begin{figure}
    \centering
    \includegraphics[width=0.5\textwidth]{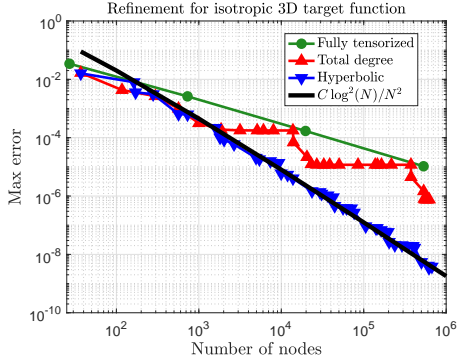}
    \caption{Isotropic refinement for trigonometric interpolation of $f_{(1,1,1)}(\bm x) = \prod_{k=1}^3 h_1(x_k)$ with different choices of $\Lambda^{\bm \alpha}(L)$. Here, $\bm \alpha=\bm 1$ and refinement occurs solely by incrementing $L$. As expected, the hyperbolic cross-section refinement converges at the expected rate and outperforms the total-degree and fully tensorized methods.}
    \label{fig:isotropic_example}
\end{figure}

The domain of interpolation for~\cref{equ:multidim_periodic_poly} is $\Gamma = [-1,1]^d$. The Fourier coefficients obey the estimate~\cref{equ:fcoef_bound_multidim}, so a hyperbolic function space like \cref{equ:hyperbolic} is appropriate, as~\cref{fig:isotropic_example} demonstrates. We calculate the error by drawing 2000 validation points $\bm x_j \sim \mathcal{U}(\Gamma)$, where $\mathcal{U}(\Gamma)$ is the uniform distribution on $\Gamma$, with
$$
\text{error} = \max_{1 \leq j \leq 2000} |f(\bm x_j)-I_{\Theta}[f](\bm x_j)| \, .
$$
In the isotropic example, the initial grids have approximately the same number of nodes, and we refine up to a maximum of 700000 nodes. 

We now consider target functions with various numbers of inputs and anisotropy. The initial grid for each refinement strategy has approximately the same number of nodes, and we refine up to a maximum of 200000 nodes.

\begin{figure}
    \centering 
	\includegraphics[width=0.5\textwidth]{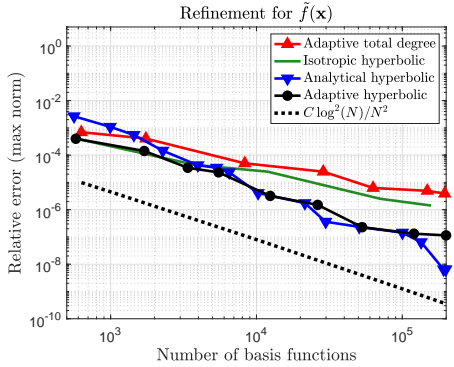}
	\caption{Convergence results for \cref{equ:6d_anisotropic_poly}. The adaptive hyperbolic cross section methods matches the convergence rate of the analytic anisotropy, but without using any prior knowledge.}
	\label{fig:polynomial_results}
\end{figure}

\begin{table}
\centering
\begin{tabular}{c|c|c|c}
	$f_{\bm i}$ & Final $\hat{\alpha}_1 / \hat{\alpha}_2$ (hyperbolic) & Final $\hat{\alpha}_1 / \hat{\alpha}_2$ (TD) & True $\alpha_1 / \alpha_2$ \\ \hline
	$(1,2)$ & 0.72 & 0.73 & 0.75 \\
	$(1,3)$ & 0.61 & 0.61 & 0.60 \\
	$(1,4)$ & 0.49 & 0.49 & 0.50 \\
	$(1,5)$ & 0.45 & 0.45 & 0.43 \\
	$(2,3)$ & 0.84 & 0.84 & 0.80 \\
	$(2,4)$ & 0.68 & 0.68 & 0.67 \\
	$(2,5)$ & 0.62 & 0.62 & 0.57 \\
	$(3,4)$ & 0.81 & 0.81 & 0.83 \\ 
	$(3,5)$ & 0.74 & 0.74 & 0.71 \\
	$(4,5)$ & 0.91 & 0.91 & 0.86 \\ \hline 
\end{tabular} 
\caption{Anisotropy ratios for two-dimensional product functions at the end of refinement. Column 2 uses~\cref{alg:adaptive_loop} and Column 3 uses~\cref{rem:total_degree}.}
\label{tbl:anisotropy_ratios}
\end{table}

Next we consider an anisotropic example. In \cref{fig:polynomial_results}, we compare different anisotropic grids, and we use the six-dimensional target function
\begin{equation}
\tilde{f}(\bm x) = h_1(x_1)h_5(x_4) + h_2(x_2) h_5(x_5) + h_3(x_3)h_5(x_6) \, .
\label{equ:6d_anisotropic_poly}
\end{equation}
Since $h_k \in H^k(\Gamma)$, then by~\cref{equ:fcoef_bound_multidim} and~\cref{equ:hyperbolic}, we know the anisotropy of $\tilde{f}$ beforehand:
$$
\bm \alpha = (1,2,3,5,5,5) + \bm 2 = (3,4,5,7,7,7).
$$
The line in~\cref{fig:polynomial_results} labeled ``Analytical hyperbolic'' uses the known anisotropy $\bm \alpha$, while the adaptive strategies solve the relevant least-squares problem for $\bm{\hat \alpha}$ at each refinement iteration. In terms of convergence behavior, all strategies with a hyperbolic cross-section space outperform the total-degree space of~\cref{rem:total_degree}. Additionally, the adaptive algorithms based on solving the least-squares problem~\cref{equ:sparse_lsq} converge at a similar rate as using the known target space $\Lambda^{\bm \alpha}(L)$ directly. Both the adaptive and analytical anisotropic strategies converge at approximately the rate given in~\cref{equ:sparse_convergence_rate}. {\em This shows that~\cref{alg:adaptive_loop} is well suited to handle periodic models where the anisotropy is not known a priori}. 

At the end of refinement, we obtain the following anisotropy estimates (normalized so that $\hat{\alpha}_1=\alpha_1=3$):
\begin{align*}
\bm{\hat \alpha}_{hyp} &= (3.00, \ 3.53, \ 4.35, \ 5.58, \ 5.70, \ 5.73),  \\
\bm{\hat \alpha}_{TD} &= (3.00, \ 3.63, \ 4.51, \ 6.11, \ 5.73, \ 5.40)\, .
\end{align*}
In~\cref{tbl:anisotropy_ratios}, we show the anisotropy ratios at the end of adaptive refinement for two-dimensional product polynomials of the form \cref{equ:multidim_periodic_poly}. We compute the true anisotropy ratio for $f_{\bm i}$ by recalling $\alpha_k = i_k+2$. Both~\cref{alg:adaptive_loop} and the modifications in~\cref{rem:total_degree} are reasonably able to detect the relative anisotropy of the target function.

\subsection{Particle in a two-dimensional box}
\label{subsec:2D_PIB}
We construct a two-dimensional particle in a box (PIB) system. This is a staple example in textbooks on quantum mechanics, e.g.~\cite{Levine_2014}. Here, the anisotropy arises from different perturbations in the $x$ and $y$ directions. As discussed in~\cite{Levine_2014}, the Hamiltonian for the unperturbed one-dimensional PIB on the interval $[0,1]$ is
\begin{equation}
\mathcal{\hat{H}} = - \frac12 \frac{\text{d}^2}{\text{d}x^2} + V(x), \qquad V(x) = \begin{cases} 0, & x \in [0,1] \\ \infty, & \text{else} \end{cases}
\label{equ:pib_unperturbed}
\end{equation}
where we have used atomic units and set the particle mass equal to the electron rest mass, $m_e=1$. For $n=1,2,\dots$, the normalized wavefunctions satisfy
$$
\mathcal{\hat{H}} \, \psi_n = E_n \, \psi_n \qquad \implies \qquad \psi_n(x) = \sqrt{2} \, \sin\left( n \pi x\right),~E_n = \frac12 {n^2 \pi^2} \, .
$$

Inspired by exercises in quantum mechanics textbooks~\cite[p.~261]{Levine_2014}, we use the potentials
$$
f_1(x) = \begin{cases} 15, & x \in [0,1/4] \cup [3/4,1] \\ 0, & \text{else} \end{cases}, \qquad f_2(y) = 60 \left(y-\frac12\right)^2
$$
in our two-dimensional perturbed PIB system and treat them as perturbations. Note that the maximum value of each perturbation is less than $E_2^{(0)}$, the energy of the unperturbed $n=2$ energy level. The two-dimensional Hamiltonian is
\begin{equation}
\mathcal{\hat H} = -\frac12 \nabla^2 + V(x) + V(y) + f_1(x) + f_2(y)
\label{equ:pib_perturbed}
\end{equation}
where $V$ is given in~\cref{equ:pib_unperturbed}. The full two-dimensional wavefunction has the form
$$
\Psi_{\bm n}(x,y) =  \psi_{n_1}(x) \, \psi_{n_2} (y)
$$
where $n_k$ is the quantum number in dimension $k$.

We will demonstrate the performance of various refinement strategies on the wavefunction of~\cref{equ:pib_perturbed} corresponding to $n_x=n_y=2$. To evaluate the target wavefunction, we first decompose~\cref{equ:pib_perturbed} into the $x$ and $y$ parts and apply first-order nondegenerate perturbation theory (see, e.g., \cite[p.~233]{Levine_2014}). The first-order correction to the wavefunction for the $x$ component is
\begin{equation}
\psi^{(1)}_{2,x}(x) = \sum_{n \neq 2} \frac{\int_0^1 \psi_n^{(0)}(u) \, f_1(u) \, \psi_2^{(0)}(u) \, \text{d}u }{E_2^{(0)}-E_n^{(0)}} \ \psi_n^{(0)}(x)
\label{equ:psi_2_x_corr}
\end{equation}
where $\psi_n^{(0)}$ and $E_n^{(0)}$ are the unperturbed wavefunctions and energies corresponding to~\cref{equ:pib_perturbed}. Similarly, for the $y$ component, we get
\begin{equation}
\psi^{(1)}_{2,y}(y) = \sum_{n \neq 2} \frac{\int_0^1 \psi_n^{(0)}(u) \, f_2(u) \, \psi_2^{(0)}(u) \, \text{d}u }{E_2^{(0)}-E_n^{(0)}} \ \psi_n^{(0)}(y) \, .
\label{equ:psi_2_y_corr}
\end{equation}
Thus, we take the two-dimensional target wavefunction as
\begin{equation}
\Psi_{2,2}(x,y) = \left(\psi_2^{(0)}(x) + \psi_{2,x}^{(1)}(x) \right) \left(\psi_2^{(0)}(y) + \psi_{2,y}^{(1)}(y) \right) \, .
\label{equ:PIB_target}
\end{equation}

We evaluate the integral coefficients in~\cref{equ:psi_2_x_corr}-\cref{equ:psi_2_y_corr} with Maple and find that the only nonzero coefficients correspond to functions of the form $\sin(2\pi k x)$, yielding an {\em a priori} anisotropy estimate. The $L^2$-Fourier coefficients of $\Psi_{2,2}(x,y)$ decay like $O(1/k^3)$ in the $x$ component and $O(1/k^5)$ in $y$, where $k$ is the coefficient index. This both justifies the use of approximation space $\Lambda^{\bm \alpha}_{hyp}$ and gives the prior anisotropy $\bm \alpha = (3,5)$. Computationally, we truncate the series in \cref{equ:psi_2_x_corr}-\cref{equ:psi_2_y_corr} at $N=10^4$ terms and use~\cref{equ:PIB_target} as our target function.  In practice, though, one would not use Fourier interpolation on a known truncated Fourier series; instead, a more accurate solution technique would provide the target wavefunction, e.g.~\cite{Consortini_Frieden_1976}. Perturbation theory, however, is straightforward enough to use for the end goal of demonstrating the convergence behavior of our adaptive refinement method.

\begin{figure}
    \centering
    \includegraphics[width=0.5\textwidth]{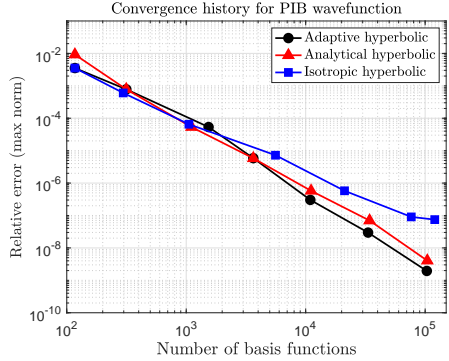}
    \caption{Convergence history of approximating~\cref{equ:PIB_target} with various techniques.}
    \label{fig:PIB_example}
\end{figure}

\begin{figure}
    \centering
    \includegraphics[width=0.49\textwidth]{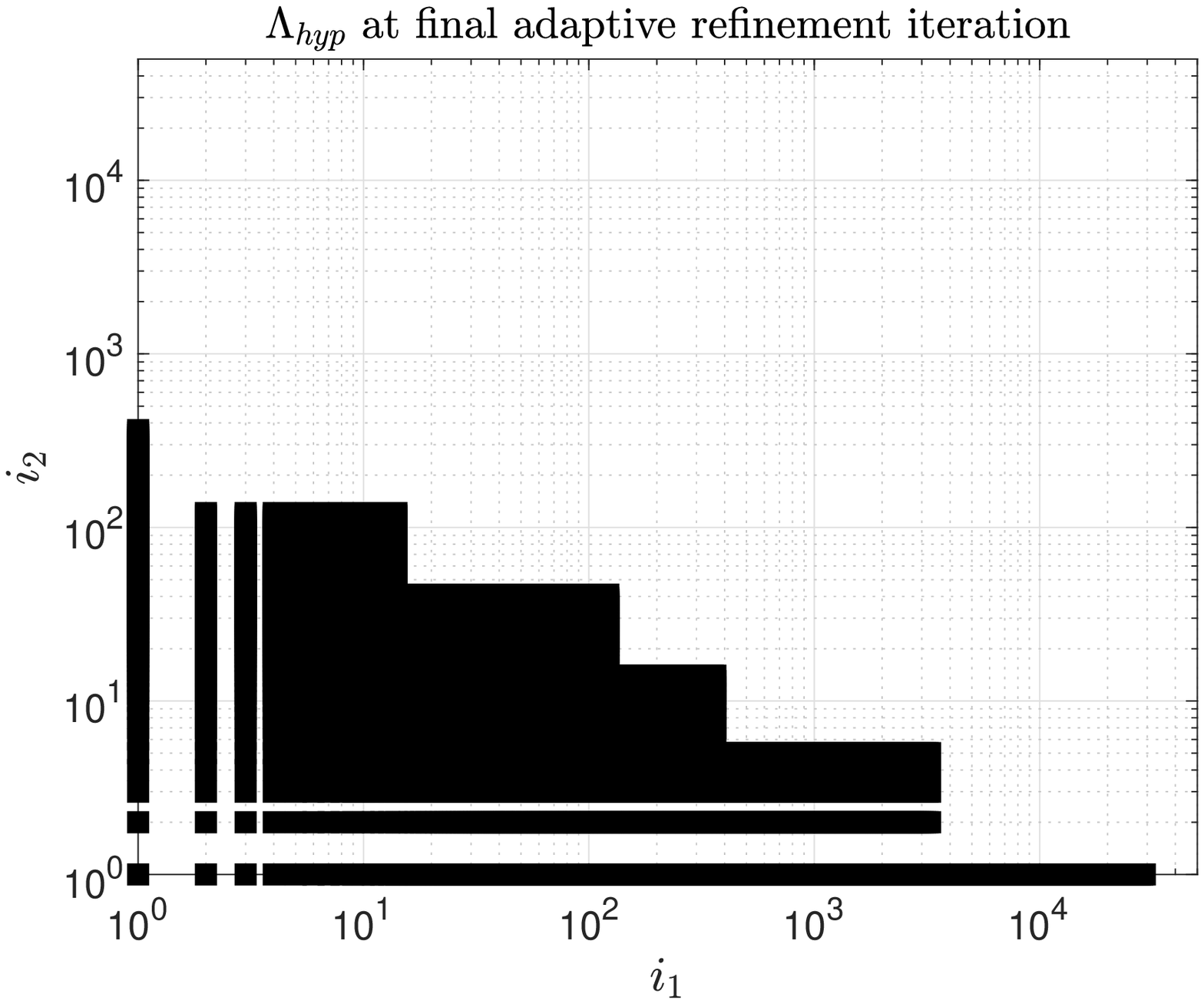}
    \includegraphics[width=0.49\textwidth]{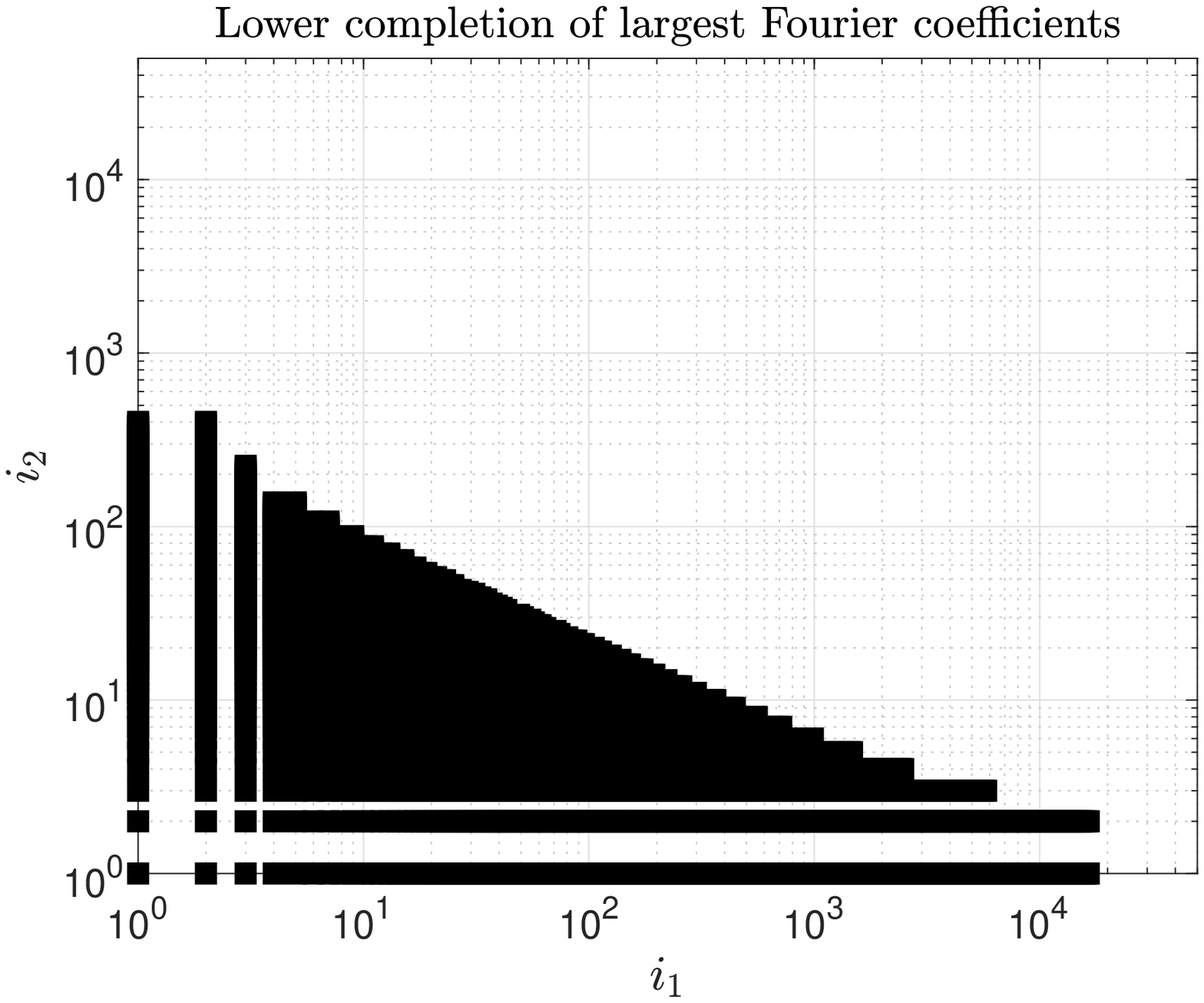}
    \caption{$\Lambda_{hyp}$ at final iteration of adaptive refinement (left); lower completion of indices for the largest Fourier coefficients (right). Note the logarithmic scaling. Equal vertical and horizontal axes are chosen to display the anisotropy.}
    \label{fig:PIB_lambda}
\end{figure}

For sparse interpolation, we use the hyperbolic index set $\Lambda_{hyp}$ and refine according to three strategies: adaptive (\cref{alg:adaptive_loop}), analytical anisotropy, and isotropic. We show the convergence behavior in~\cref{fig:PIB_example}. Similarly to~\cref{subsec:results_polynomials}, adaptive refinement performs as well as analytical anisotropic refinement, but without any prior knowledge. Asymptotically, the errors of the analytical anisotropic and adaptive strategies in~\cref{fig:PIB_example} decay at roughly the same rate and are an order of magnitude better than isotropic refinement. Furthermore, since we know the Fourier coefficients explicitly as a result of~\cref{equ:psi_2_x_corr}-\cref{equ:psi_2_y_corr}, we may construct the optimal lower approximation space directly from the explicit coefficients. In~\cref{fig:PIB_lambda} we show $\Lambda_{hyp}$ at the final iteration of adaptive refinement along with the smallest lower set of size $\#(\Lambda_{hyp})$ containing the $N \leq \#(\Lambda_{hyp})$ largest Fourier coefficients. \Cref{fig:PIB_lambda} shows that adaptive refinement closely resembles the lower set containing the $N \leq \# (\Lambda_{hyp})$ largest Fourier coefficients, except for some rectangular gaps introduced by the growth rule $m(l) = 3^l$. We chose equal vertical and horizontal axes to make the anisotropy clear.

\subsection{The 2-butene potential energy surface}
\label{subsec:chemistry}

\begin{figure}
	\centering
	\includegraphics[width=0.4\textwidth]{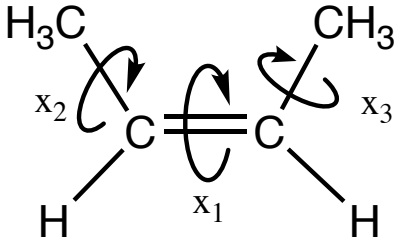}
	\caption{Molecular structure of 2-butene, labeled with rotations of interest.}
	\label{fig:2-butene}
\end{figure}

Now we consider the motivating application of this paper: the adaptive approximation of a molecule's potential energy surface (PES) where the anisotropy $\bm \alpha$ is not known beforehand. The molecule of interest is 2-butene, whose molecular structure is shown in~\cref{fig:2-butene}. 

At the quantum-mechanical level, the energy $\mathcal{E}_n$ of a molecule with an arrangement of nuclei described by $\bm q$ satisfies the Schr\"odinger equation~\cite{Levine_2014}:
\begin{equation}
\mathcal{\hat H}(\bm q) \, \Psi_n (\bm y; \bm q) = \mathcal{E}_n(\bm q) \, \Psi_n(\bm y; \bm q) \, .
\label{equ:schrodinger}
\end{equation}
Above, $\mathcal{\hat H}$ is the molecular Hamiltonian operator, $\mathcal{E}_n$ is the energy of electronic state $n \geq 0$, and $\Psi_n$ is the (possibly complex-valued) wavefunction of state $n$ as a function of electron position $\bm y$. All terms depend parametrically on the nuclear geometry $\bm q$. Physically, $|\Psi_n(\bm y; \bm q)|^2$ is the probability distribution function of observing an electron of energy state $n$ at position $\bm y$ in a molecule of geometry $\bm q$. As a function of geometry, $\mathcal{E}_n(\bm q)$ is the PES corresponding to energy state $n$.

For an $N$-atom molecule, one may express $\bm q$ in Cartesian coordinates as a vector with $3N$ components or in internal coordinates (bond lengths, bond angles, and torsion angles) as a vector of $3N-6$ components. We opt for the latter, which has fewer components and directly enables the varying of geometric features. Only a few of the geometry components, denoted by $\bm x$ (the \textit{design variables}), may be needed in a particular study; we optimize over the rest, $\bm \xi$ (the \textit{remainder variables}):
\begin{equation} 
E_n(\bm x) = \min_{\bm \xi} \, \mathcal{E}_n(\bm x, \bm \xi) \, .
\label{equ:relaxed_PES}
\end{equation} 
$E_n(\bm x)$ is called the \textit{relaxed PES} for state $n$. Chemical intuition and knowledge of the system guides the selection of design variables.  Furthermore, for rotational design variables $\bm \theta$, a polynomial interpolant does not guarantee periodicity of $\nabla E_n$ with respect to $\bm \theta$, which leads to nonphysical phenomena (e.g., nonconservation of energy). Therefore, a trigonometric interpolation basis is appropriate when $\bm x$ contains only bond angles and torsion angles.\footnote{Bond lengths, in general, are not periodic over an interpolation domain, so approximation by trigonometric polynomials would lead to inaccuracies at the domain boundary~\cite{Helmberg_1994}.}

As hinted earlier, solving the optimization~\cref{equ:relaxed_PES} subject to the generalized eigenvalue problem~\cref{equ:schrodinger} is a prohibitively expensive calculation. To trim down computational cost, it is common practice in quantum chemistry to use approximate Hamiltonians and wavefunctions~\cite{Levine_2014}. In our case, we use density functional theory (with the B3LYP hybrid functional) to simplify the Hamiltonian~\cite{Hohenberg_Kohn_1964, Kohn_Sham_1965,Stephens_et_al_1994}, and we approximate the wavefunctions with the 6-311G* Pople basis set~\cite{Krishnan_et_al_1980}. We use the Gaussian 16 software package~\cite{Gaussian_16} to handle the approximation of Hamiltonians and wavefunctions. By default, Gaussian 16 performs the optimization in~\cref{equ:relaxed_PES} using a variant of the EDIIS algorithm tuned for molecular geometry optimizations~\cite{Li_Frisch_2006}.

Previous work constructed a sparse polynomial interpolant of $E_0(\bm x)$ and $E_1(\bm x)$ for 2-butene to study the transition from the \textit{cis}- to \textit{trans}- conformation via the first singlet excited state~\cite{Nance_Jakubikova_Kelley_2014}. We use the same design variables from that study, shown in~\cref{fig:2-butene}. The design variable $x_1$ is more influential on the PES than $x_2$ and $x_3$, but the exact anisotropy is not known in advance.

\begin{figure}
    \centering
    \includegraphics[width=0.5\textwidth]{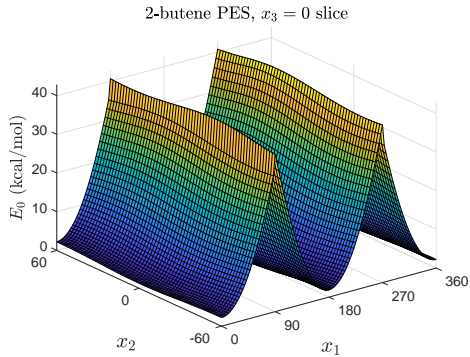}
    \caption{Slice of 2-butene PES for $x_3=0$.}
    \label{fig:2-butene_PES}
\end{figure}

The domain for our 2-butene ground-state ($n=0$) PES is $\Gamma = [0,360] \times [-60,60] \times [-60,60]$. The coordinates $x_2$ and $x_3$ correspond to dihedral rotations of $\text{CH}_3$, which have period $120^\circ$. We show a slice of the 2-butene PES in \cref{fig:2-butene_PES}. The ridges at $x_1=90$ and $x_1=270$ indicate a discontinuous first derivative, so we hypothesize that the hyperbolic function space \cref{equ:hyperbolic} and \cref{alg:adaptive_loop} are appropriate for this problem. 

There are numerous sources of noise going into the evaluation of $E_0(\bm x)$: density functional theory approximates the Hamiltonian, the 6-311G* basis set approximates the wavefunction $\Psi$, and the optimization~\cref{equ:relaxed_PES} has internal stopping criteria. Therefore, we do not report the max error of the interpolant $I_{\Theta}[E_0](\bm x)$, which could be heavily skewed by non-interpolatory error. Instead, we give the root-mean-square error (RMSE) over 2000 validation points drawn uniformly over $\Gamma$:
\begin{equation} 
\text{RMSE} = \sqrt{ \frac{\sum_{j=1}^{2000} (E_0(\bm x_j)-I_{\Theta}[E_0](\bm x_j))^2}{2000}}\, , \hspace{2em} \bm x_j \sim \mathcal{U}(\Gamma)\, .
\label{equ:rmse}
\end{equation}

\begin{figure}
\centering
\includegraphics[width=0.49\textwidth]{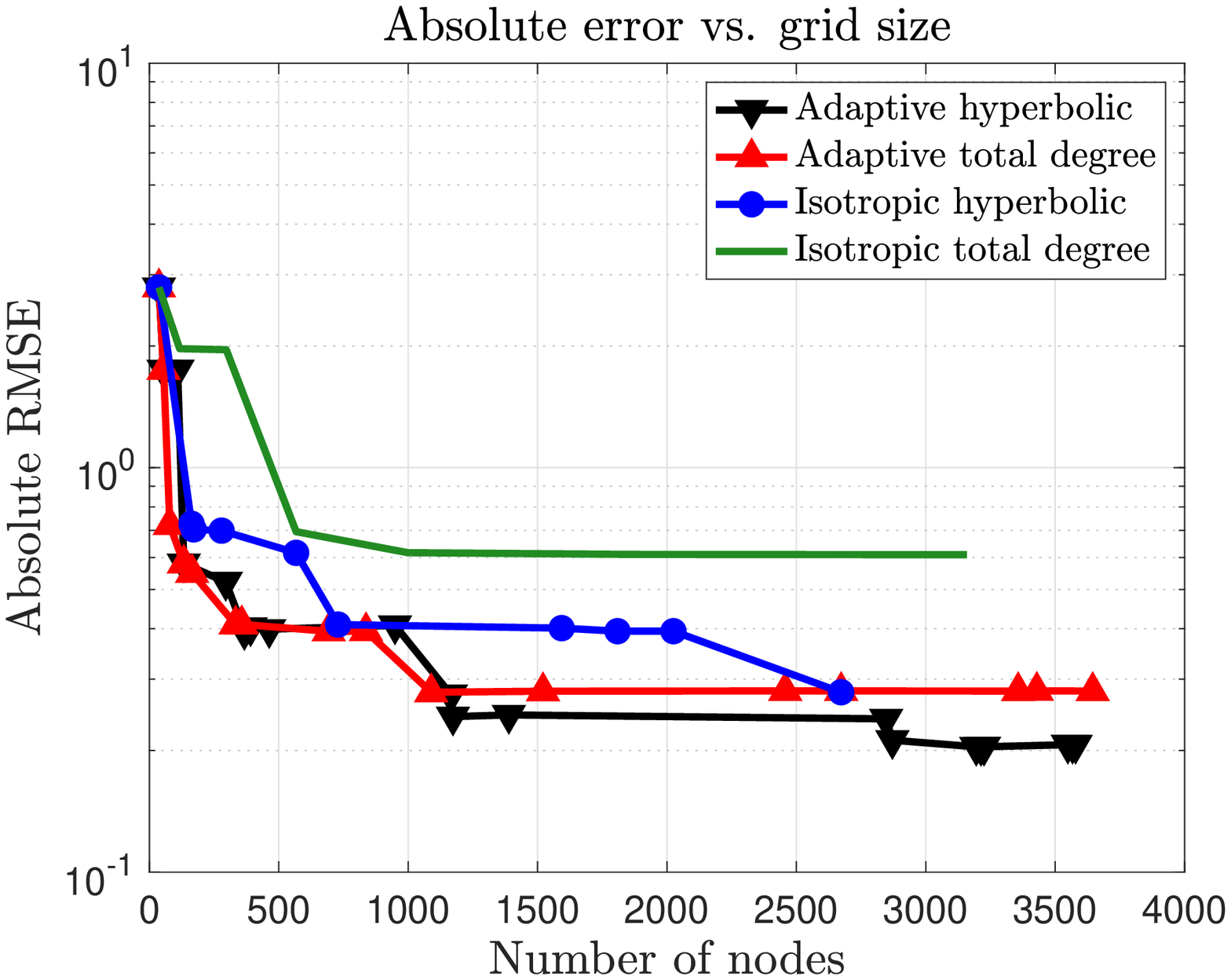}
\includegraphics[width=0.49\textwidth]{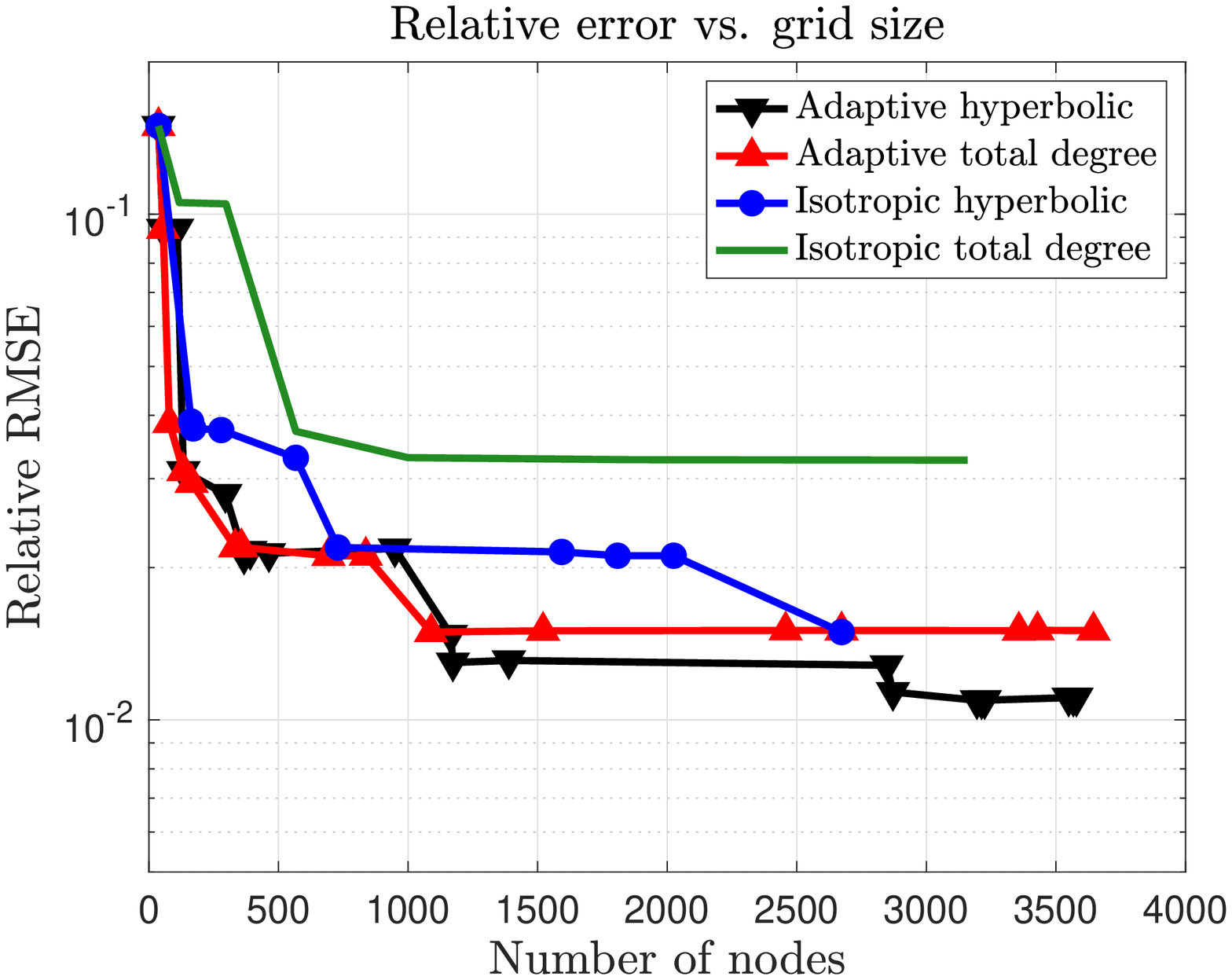}
\caption{Absolute (left) and relative (right) error results for sparse interpolation of 2-butene ground-state PES.}
\label{fig:PES_results}
\end{figure}

In our sparse grid constructions, we use $\Theta_{opt}$ based on the hyperbolic function space $\Lambda$ in~\cref{equ:hyperbolic} as well as the total-degree space~\cref{equ:total_degree}. For each variety of $\Theta$, we refine both adaptively (according to~\cref{alg:adaptive_loop} or~\cref{rem:total_degree}) and isotropically (taking $\bm \alpha=\bm 1$ and incrementing $L$). In all cases, we initialize each grid with 37 nodes. Each function sample takes approximately 30 seconds to evaluate, and occasionally the optimization~\cref{equ:relaxed_PES} may fail to converge to the correct (or any) local minimum. Due to limitations on available computing time, we refine up to a maximum of only 4000 nodes. If a refinement strategy terminates prior to 4000 nodes, that is because the next increment of $L$ would result in the number of nodes exceeding 4000. We show the results in~\cref{fig:PES_results}. Following Pople in his 1998 Nobel lecture, we adopt 1 kcal/mol as the threshold of acceptable chemical accuracy for energies~\cite{Pople_1999}.

First, we note that the asymptotic absolute RMS errors in~\cref{fig:PES_results} are consistent with Pople's definition of chemical accuracy for energies (i.e., less than 1 kcal/mol). Furthermore, the limiting relative error for adaptive hyperbolic refinement is approximately 1\%. Second, even though the one-dimensional interpolation rule~\cref{equ:rule_1d} grows exponentially, we can still add smaller batches of nodes at each iteration by using~\cref{equ:theta_opt} and~\cref{alg:adaptive_loop}, which mitigates the exponential growth of \cref{equ:rule_1d}. Third, in both the adaptive and isotropic cases, the asymptotic error is lower for a hyperbolic cross-section than for a total-degree space.

\section{Conclusion and future work}
\label{sec:conclusion}

In this work, we have presented a quasi-optimal dimensionally adaptive method for sparse interpolation with a trigonometric basis. Our approach targets applications where the surrogate models must be periodicity-preserving and where the anisotropy is not known beforehand. For target functions of known finite smoothness, our algorithm matches the theoretical convergence rate, outperforms the total-degree space asymptotically, and produces a good approximation to the anisotropy. The open-source and freely available Tasmanian package contains a user-friendly implementation.

In the future, we will apply adaptive refinement to more complicated chemical systems. In particular, we aim to approximate potential energy surfaces where the geometry domain includes bond lengths, bond angles, and dihedral angles---a mix of periodic and nonperiodic inputs. To do this, we will apply trigonometric and polynomial interpolation to the periodic and nonperiodic components, respectively.

\section*{Acknowledgements}
\label{sec:ack}
The first author was supported by an NSF Graduate Research Fellowship under DGE-1746939 and by an appointment to the Oak Ridge National Laboratory Advanced Short-Term Research Opportunity (ASTRO) Program, sponsored by the U.S. Department of Energy and administered by the Oak Ridge Institute for Science and Education. The second author was supported by the Exascale Computing Project (17-SC-20-SC), a collaborative effort of the U.S. Department of Energy Office of Science and the National Nuclear Security Administration; 
the U.S. Defense Advanced Research Projects Agency, Defense Sciences Office under contract and award numbers HR0011619523 and 1868-A017-15;
and by the Laboratory Directed Research and Development program at the Oak Ridge National Laboratory, which is operated by UT-Battelle, LLC, for the U.S. Department of Energy under Contract DE-AC05-00OR22725.

The authors also acknowledge the use of the High Performance Computing Center at North Carolina State University.


\bibliographystyle{siamplain}
\bibliography{main}

\begin{thebibliography}{10}

\bibitem{Barthelmann_Novak_Ritter_2000}
{\sc V.~Barthelmann, E.~Novak, and K.~Ritter}, {\em High dimensional polynomial
  interpolation on sparse grids}, Adv. Comput. Math., 12 (2000), pp.~273--288,
  \url{https://doi.org/10.1023/A:1018977404843}.

\bibitem{Beck_et_al_2014}
{\sc J.~Beck, F.~Nobile, L.~Tamellini, and R.~Tempone}, {\em Convergence of
  quasi-optimal stochastic {G}alerkin methods for a class of {PDE}s with random
  coefficients}, Comput. Math. Appl., 67 (2014), pp.~732--751,
  \url{https://doi.org/10.1016/j.camwa.2013.03.004}.

\bibitem{Bellman_1961}
{\sc R.~Bellman}, {\em Adaptive Control Processes: A Guided Tour}, Princeton
  University Press, 1961.

\bibitem{Bungartz_Griebel_2004}
{\sc H.-J. Bungartz and M.~Griebel}, {\em Sparse grids}, Acta Numer., 13
  (2004), pp.~147--269, \url{https://doi.org/10.1017/S0962492904000182}.

\bibitem{Chkifa_Cohen_Schwab_2014}
{\sc A.~Chkifa, A.~Cohen, and C.~Schwab}, {\em High-dimensional adaptive sparse
  polynomial interpolation and applications to parametric {PDE}s}, Found.
  Comput. Math., 14 (2014), pp.~601--633,
  \url{https://doi.org/10.1007/s10208-013-9154-z}.

\bibitem{Clenshaw_Curtis_1960}
{\sc C.~W. Clenshaw and A.~R. Curtis}, {\em A method for numerical integration
  on an automatic computer}, Numer. Math., 2 (1960), pp.~197--205,
  \url{https://doi.org/10.1007/BF01386223}.

\bibitem{Consortini_Frieden_1976}
{\sc A.~Consortini and B.~Frieden}, {\em Quantum-mechanical solution for the
  simple harmonic oscillator in a box}, Nuov Cim B, 35 (1976), p.~153–164,
  \url{https://doi.org/10.1007/BF02724052}.

\bibitem{Dzyadyk_Dzyadyk_Prypik_1981}
{\sc V.~K. Dzyadyk, {S. Y. Dzyadyk}, and A.~S. Prypik}, {\em Asymptotic
  behavior of {L}ebesgue constants in trigonometric interpolation}, Ukrainian
  Math. J., 33 (1981), pp.~553--559, \url{https://doi.org/10.1007/BF01085428}.

\bibitem{Ehlich_Zeller_1966}
{\sc H.~Ehlich and K.~Zeller}, {\em {A}uswertung der {N}ormen von
  {I}nterpolationsoperatoren}, Math. Ann., 164 (1966), pp.~105--112,
  \url{https://doi.org/10.1007/BF01429047}.

\bibitem{EldredBurkardt2009}
{\sc M.~Eldred and J.~Burkardt}, {\em Comparison of non-intrusive polynomial
  chaos and stochastic collocation methods for uncertainty quantification},
  AIAA, 2009, pp.~1--20, \url{https://doi.org/10.2514/6.2009-976}.

\bibitem{EldredWebsterConstantine2008}
{\sc M.~Eldred, C.~Webster, and P.~Constantine}, {\em Evaluation of
  non-intrusive approaches for Wiener-Askey generalized polynomial chaos},
  AIAA, 2008, ch.~2008-1892, pp.~1--22,
  \url{https://doi.org/10.2514/6.2008-1892}.

\bibitem{Gaussian_16}
{\sc M.~J. Frisch, G.~W. Trucks, H.~B. Schlegel, G.~E. Scuseria, M.~A. Robb,
  J.~R. Cheeseman, G.~Scalmani, V.~Barone, G.~A. Petersson, H.~Nakatsuji,
  X.~Li, M.~Caricato, A.~V. Marenich, J.~Bloino, B.~G. Janesko, R.~Gomperts,
  B.~Mennucci, H.~P. Hratchian, J.~V. Ortiz, A.~F. Izmaylov, J.~L. Sonnenberg,
  D.~Williams-Young, F.~Ding, F.~Lipparini, F.~Egidi, J.~Goings, B.~Peng,
  A.~Petrone, T.~Henderson, D.~Ranasinghe, V.~G. Zakrzewski, J.~Gao, N.~Rega,
  G.~Zheng, W.~Liang, M.~Hada, M.~Ehara, K.~Toyota, R.~Fukuda, J.~Hasegawa,
  M.~Ishida, T.~Nakajima, Y.~Honda, O.~Kitao, H.~Nakai, T.~Vreven,
  K.~Throssell, J.~A. Montgomery, {Jr.}, J.~E. Peralta, F.~Ogliaro, M.~J.
  Bearpark, J.~J. Heyd, E.~N. Brothers, K.~N. Kudin, V.~N. Staroverov, T.~A.
  Keith, R.~Kobayashi, J.~Normand, K.~Raghavachari, A.~P. Rendell, J.~C.
  Burant, S.~S. Iyengar, J.~Tomasi, M.~Cossi, J.~M. Millam, M.~Klene, C.~Adamo,
  R.~Cammi, J.~W. Ochterski, R.~L. Martin, K.~Morokuma, O.~Farkas, J.~B.
  Foresman, and D.~J. Fox}, {\em Gaussian 16 {R}evision {A}.03}, 2016.
\newblock Gaussian Inc. Wallingford, CT.

\bibitem{Gautschi_2012}
{\sc W.~Gautschi}, {\em Numerical Analysis}, Birkh\"{a}user Basel, 2012.

\bibitem{Grafakos_2014}
{\sc L.~Grafakos}, {\em Classical Fourier Analysis}, Springer, 2014.

\bibitem{Griebel_Hamaekers_2014}
{\sc M.~Griebel and J.~Hamaekers}, {\em Fast discrete {F}ourier transform on
  generalized sparse grids}, in Sparse Grids and Applications -- Munich 2012,
  J.~Garcke and D.~Pfl\"{u}ger, eds., vol.~97 of Lecture Notes in Computational
  Science and Engineering, Springer International, 2014, pp.~75--107,
  \url{https://doi.org/10.1007/978-3-319-04537-5_4}.

\bibitem{gunzburger2014stochastic}
{\sc M.~D. Gunzburger, C.~G. Webster, and G.~Zhang}, {\em Stochastic finite
  element methods for partial differential equations with random input data},
  Acta Numer., 23 (2014), pp.~521--650,
  \url{https://doi.org/10.1017/S0962492914000075}.

\bibitem{Hallatschek_1992}
{\sc K.~Hallatschek}, {\em Fouriertransformation auf d{\"u}nnen {G}ittern mit
  hierarchischen {B}asen}, Numer. Math., 63 (1992), pp.~83--97,
  \url{https://doi.org/10.1007/BF01385849}.

\bibitem{Hart_Alexanderian_Gremaud_2017}
{\sc J.~Hart, A.~Alexanderian, and P.~Gremaud}, {\em Efficient computation of
  {S}obol' indices for stochastic models}, SIAM J. Sci. Comput., 39 (2017),
  pp.~A1514--A1530, \url{https://doi.org/10.1137/16M106193X}.

\bibitem{Helmberg_1994}
{\sc G.~Helmberg}, {\em The {G}ibbs phenomenon for {F}ourier interpolation}, J.
  Approx. Theory, 78 (1994), pp.~41--63,
  \url{https://doi.org/10.1006/jath.1994.1059}.

\bibitem{Hohenberg_Kohn_1964}
{\sc P.~Hohenberg and W.~Kohn}, {\em Inhomogeneous electron gas}, Phys. Rev.,
  136 (1964), pp.~B864--B871, \url{https://doi.org/10.1103/PhysRev.136.B864}.

\bibitem{Jackson_1930}
{\sc D.~Jackson}, {\em The Theory of Approximation}, American Mathematical
  Society, 1930.

\bibitem{jakeman2011characterization}
{\sc J.~D. Jakeman, R.~Archibald, and D.~Xiu}, {\em Characterization of
  discontinuities in high-dimensional stochastic problems on adaptive sparse
  grids}, J. Comput. Phys., 230 (2011), pp.~3977--3997,
  \url{https://doi.org/10.1016/j.jcp.2011.02.022}.

\bibitem{jakeman2013minimal}
{\sc J.~D. Jakeman, A.~Narayan, and D.~Xiu}, {\em Minimal multi-element
  stochastic collocation for uncertainty quantification of discontinuous
  functions}, J. Comput. Phys., 242 (2013), pp.~790--808,
  \url{https://doi.org/10.1016/j.jcp.2013.02.035}.

\bibitem{jakeman2012local}
{\sc J.~D. Jakeman and S.~G. Roberts}, {\em Local and dimension adaptive
  stochastic collocation for uncertainty quantification}, in Sparse Grids and
  Applications, Springer, 2012, pp.~181--203,
  \url{https://doi.org/10.1007/978-3-642-31703-3_9}.

\bibitem{Katznelson_2004}
{\sc Y.~Katznelson}, {\em An Introduction to Harmonic Analysis}, Cambridge
  University Press, 3rd~ed., 2004.

\bibitem{khakhutskyy2016spatially}
{\sc V.~Khakhutskyy and M.~Hegland}, {\em Spatially-dimension-adaptive sparse
  grids for online learning}, in Sparse Grids and Applications -- Stuttgart
  2014, Springer, 2016, pp.~133--162,
  \url{https://doi.org/10.1007/978-3-319-28262-6_6}.

\bibitem{klimke2005algorithm}
{\sc A.~Klimke and B.~Wohlmuth}, {\em Algorithm 847: spinterp: {P}iecewise
  multilinear hierarchical sparse grid interpolation in {MATLAB}}, ACM Trans.
  Math. Software, 31 (2005), pp.~561--579,
  \url{https://doi.org/10.1145/1114268.1114275}.

\bibitem{Kohn_Sham_1965}
{\sc W.~Kohn and L.~J. Sham}, {\em Self-consistent equations including exchange
  and correlation effects}, Phys. Rev., 140 (1965), pp.~A1133--A1138,
  \url{https://doi.org/10.1103/PhysRev.140.A1133}.

\bibitem{Kreyszig_1978}
{\sc E.~Kreyszig}, {\em Introductory Functional Analysis with Applications},
  John Wiley and Sons, 1978.

\bibitem{Krishnan_et_al_1980}
{\sc R.~Krishnan, J.~S. Binkley, R.~Seeger, and J.~A. Pople}, {\em
  Self‐consistent molecular orbital methods. {XX}. {A} basis set for
  correlated wave functions}, J. Chem. Phys., 72 (1980), pp.~650--654,
  \url{https://doi.org/10.1063/1.438955}.

\bibitem{Kritzer_Pillichshammer_Wozniakowski_2014}
{\sc P.~Kritzer, F.~Pillichshammer, and H.~Woźniakowski}, {\em Multivariate
  integration of infinitely many times differentiable functions in weighted
  {K}orobov spaces}, Math. Comp., 83 (2014), pp.~1189--1206,
  \url{https://doi.org/10.1090/S0025-5718-2013-02739-1}.

\bibitem{Levine_2014}
{\sc I.~N. Levine}, {\em Quantum Chemistry}, Pearson, 7th~ed., 2014.

\bibitem{Li_Frisch_2006}
{\sc X.~Li and M.~J. Frisch}, {\em Energy-represented direct inversion in the
  iterative subspace within a hybrid geometry optimization method}, J. Chem.
  Theory Comput., 2 (2006), pp.~835--839,
  \url{https://doi.org/10.1021/ct050275a}.

\bibitem{ma2009adaptive}
{\sc X.~Ma and N.~Zabaras}, {\em An adaptive hierarchical sparse grid
  collocation algorithm for the solution of stochastic differential equations},
  J. Comput. Phys., 228 (2009), pp.~3084--3113,
  \url{https://doi.org/10.1016/j.jcp.2009.01.006}.

\bibitem{Nance_Jakubikova_Kelley_2014}
{\sc J.~Nance, E.~Jakubikova, and C.~T. Kelley}, {\em Reaction path following
  with sparse interpolation}, J. Chem. Theory Comput., 10 (2014),
  pp.~2942--2949, \url{https://doi.org/10.1021/ct5004669}.

\bibitem{narayan2014adaptive}
{\sc A.~Narayan and J.~D. Jakeman}, {\em Adaptive {L}eja sparse grid
  constructions for stochastic collocation and high-dimensional approximation},
  SIAM J. Sci. Comput., 36 (2014), pp.~A2952--A2983,
  \url{https://doi.org/10.1137/140966368}.

\bibitem{Nobile_Tamellini_Tempone_2016}
{\sc F.~Nobile, L.~Tamellini, and R.~Tempone}, {\em Convergence of
  quasi-optimal sparse-grid approximation of {H}ilbert-space-valued functions:
  application to random elliptic {PDE}s}, Numer. Math., 134 (2016),
  pp.~343--388, \url{https://doi.org/10.1007/s00211-015-0773-y}.

\bibitem{nobile2008anisotropic}
{\sc F.~Nobile, R.~Tempone, and C.~G. Webster}, {\em An anisotropic sparse grid
  stochastic collocation method for partial differential equations with random
  input data}, SIAM J. Numer. Anal., 46 (2008), pp.~2411--2442,
  \url{https://doi.org/10.1137/070680540}.

\bibitem{nobile2008sparse}
{\sc F.~Nobile, R.~Tempone, and C.~G. Webster}, {\em A sparse grid stochastic
  collocation method for partial differential equations with random input
  data}, SIAM J. Numer. Anal., 46 (2008), pp.~2309--2345,
  \url{https://doi.org/10.1137/060663660}.

\bibitem{Novak_Ritter_1996}
{\sc E.~Novak and K.~Ritter}, {\em High dimensional integration of smooth
  functions over cubes}, Numer. Math., 75 (1996), pp.~79--97,
  \url{https://doi.org/10.1007/s002110050231}.

\bibitem{Novak_Ritter_1999}
{\sc E.~Novak and K.~Ritter}, {\em Simple cubature formulas with high
  polynomial exactness}, Constr. Approx., 15 (1999), pp.~499--522,
  \url{https://doi.org/10.1007/s003659900119}.

\bibitem{Novak_Wozniakowski_2008}
{\sc E.~Novak and H.~Woźniakowski}, {\em Tractability of Multivariate
  Problems}, European Mathematical Society, 2008.

\bibitem{Papageorgiou_Wozniakowski_2010}
{\sc A.~Papageorgiou and H.~Woźniakowski}, {\em Tractability through
  increasing smoothness}, J. Complexity, 26 (2010), pp.~409--421,
  \url{https://doi.org/10.1016/j.jco.2009.12.004}.

\bibitem{pflueger12spatially}
{\sc D.~Pfl{\"u}ger}, {\em Spatially adaptive refinement}, in Sparse Grids and
  Applications, Springer, 2012, pp.~243--262,
  \url{https://doi.org/10.1007/978-3-642-31703-3_12}.

\bibitem{pfluger2010spatially}
{\sc D.~Pfl{\"u}ger, B.~Peherstorfer, and H.-J. Bungartz}, {\em Spatially
  adaptive sparse grids for high-dimensional data-driven problems}, J.
  Complexity, 26 (2010), pp.~508--522,
  \url{https://doi.org/10.1016/j.jco.2010.04.001}.

\bibitem{Pinkus_2003}
{\sc A.~Pinkus}, {\em Negative theorems in approximation theory}, Amer. Math.
  Monthly, 110 (2003), pp.~900--911,
  \url{https://doi.org/10.1080/00029890.2003.11920030}.

\bibitem{Pople_1999}
{\sc J.~A. Pople}, {\em {N}obel lecture: {Q}uantum chemical models}, Rev. Mod.
  Phys., 71 (1999), pp.~1267--1274,
  \url{https://doi.org/10.1103/RevModPhys.71.1267}.

\bibitem{Rivlin_1974}
{\sc T.~J. Rivlin}, {\em The Chebyshev Polynomials}, Wiley, 1st~ed., 1974.

\bibitem{Sickel_Ullrich_2009}
{\sc W.~Sickel and T.~Ullrich}, {\em Tensor products of {S}obolev--{B}esov
  spaces and applications to approximation from the hyperbolic cross}, J.
  Approx. Theory, 161 (2009), pp.~748--786,
  \url{https://doi.org/10.1016/j.jat.2009.01.001}.

\bibitem{Smolyak_1963}
{\sc S.~A. Smolyak}, {\em Quadrature and interpolation formulas for tensor
  products of certain classes of functions}, Dokl. Akad. Nauk SSSR, 148 (1963),
  pp.~1042--1045.

\bibitem{Stephens_et_al_1994}
{\sc P.~J. Stephens, F.~J. Devlin, C.~F. Chabalowski, and M.~J. Frisch}, {\em
  Ab initio calculation of vibrational absorption and circular dichroism
  spectra using density functional force fields}, J. Phys. Chem., 98 (1994),
  pp.~11623--11627, \url{https://doi.org/10.1021/j100096a001}.

\bibitem{Stoer_Bulirsch_1993}
{\sc J.~Stoer and R.~Bulirsch}, {\em Introduction to Numerical Analysis},
  Springer-Verlag, 2nd~ed., 1993.
\newblock Translated from German by R. Bartels, W. Gautschi, and C. Witzgall.

\bibitem{stoyanov2018adaptive}
{\sc M.~Stoyanov}, {\em Adaptive sparse grid construction in a context of local
  anisotropy and multiple hierarchical parents}, in Sparse Grids and
  Applications -- Miami 2016, Springer, 2018, pp.~175--199,
  \url{https://doi.org/10.1007/978-3-319-75426-0_8}.

\bibitem{Tasmanian}
{\sc M.~Stoyanov}, {\em {TASMANIAN} sparse grids (version 6.0)}, Tech. Report
  ORNL/TM-2015/596, Oak Ridge National Laboratory, 2018,
  \url{https://tasmanian.ornl.gov/}.

\bibitem{stoyanov2017predicting}
{\sc M.~Stoyanov, P.~Seleson, and C.~Webster}, {\em Predicting fracture
  patterns in simulations of brittle materials under variable load and material
  strength}, in 19th AIAA Non-Deterministic Approaches Conference, 2017,
  p.~1326, \url{https://doi.org/10.2514/6.2017-1326}.

\bibitem{Stoyanov_Webster_2015}
{\sc M.~Stoyanov and C.~G. Webster}, {\em A gradient-based sampling approach
  for dimension reduction of partial differential equations with stochastic
  coefficients}, Int. J. Uncertain. Quantif., 5 (2015), pp.~49--72,
  \url{https://doi.org/10.1615/Int.J.UncertaintyQuantification.2014010945}.

\bibitem{Stoyanov_Webster_2016}
{\sc M.~K. Stoyanov and C.~G. Webster}, {\em A dynamically adaptive sparse
  grids method for quasi-optimal interpolation of multidimensional functions},
  Comput. Math. Appl., 71 (2016), pp.~2449--2465,
  \url{https://doi.org/10.1016/j.camwa.2015.12.045}.

\bibitem{Todor_Schwab_2007}
{\sc R.~A. Todor and C.~Schwab}, {\em {Convergence rates for sparse chaos
  approximations of elliptic problems with stochastic coefficients}}, IMA J.
  Numer. Anal., 27 (2007), pp.~232--261,
  \url{https://doi.org/10.1093/imanum/drl025}.

\bibitem{Tran_Webster_Zhang_2017}
{\sc H.~Tran, C.~G. Webster, and G.~Zhang}, {\em Analysis of quasi-optimal
  polynomial approximations for parameterized {PDE}s with deterministic and
  stochastic coefficients}, Numer. Math., 137 (2017), pp.~451--493,
  \url{https://doi.org/10.1007/s00211-017-0878-6}.

\end{thebibliography}

\end{document}